\documentclass{amsart}

\usepackage{amscd,amssymb}

\usepackage{hyperref}

\def\Z{{\mathbb Z}}
\def\Q{{\mathbb Q}}

\def\C{{\mathbb C}}
\def\P{{\mathbb P}}

\def\A{{\mathcal A}}

\def\K{{\mathcal K}}

\def\O{{\mathcal O}}

\def\U{{\mathcal U}}
\def\cD{{\mathcal D}}
\def\cG{{\mathcal G}}
\def\cP{{\mathcal P}}

\def\X{{\mathcal X}}

\def\cR{{\mathcal R}}

\def\a{{\mathfrak a}}
\def\f{{\mathfrak f}}
\def\g{{\mathfrak g}}

\def\k{{\mathfrak k}}
\def\n{{\mathfrak n}}
\def\p{{\mathfrak p}}
\def\r{{\mathfrak r}}

\def\u{{\mathfrak u}}

\def\G{{\Gamma}}

\def\Ghat{{\hat{\G}}}

\def\cGtilde{\tilde{\cG}}
\def\Rbar{{\overline R}}

\def\Qbar{{\overline{\Q}}}
\def\Fbar{\overline{F}}
\def\Xbar{\overline{X}}

\def\rhotilde{{\tilde{\rho}}}
\def\rhohat{{\hat{\rho}}}
\def\phihat{{\hat{\phi}}}
\def\Phihat{{\hat{\Phi}}}

\def\ntilde{{\tilde{\n}}}
\def\thetatilde{{\tilde{\theta}}}
\def\cGhat{{\widehat{\cG}}}
\def\Ntilde{{\tilde{N}}}

\def\Uhat{{\widehat{\U}}}
\def\wtilde{{\tilde{w}}}

\def\l{{\ell}}

\def\un{\mathrm{un}}

\def\prol{{(\ell)}}

\def\ur{{\mathrm{ur}}}

\def\Ql{{\Q_\ell}}
\def\Qlx{{\Q_\ell^\times}}
\def\Zl{{\Z_\ell}}
\def\Zlx{{\Z_\ell^\times}}
\def\Gl{{G_\ell}}

\def\Al{{\A_\ell}}
\def\Kl{{\K_\ell}}
\def\al{{\a_\ell}}
\def\kl{{\k_\ell}}

\def\AFS{{\A_{F,S}}}

\def\KFS{{\K_{F,S}}}
\def\aFS{{\a_{F,S}}}
\def\kFS{{\k_{F,S}}}

\def\GFS{{G_{F,S}}}
\def\OFS{{\O_{F,S}}}

\def\GFSl{{G_{F,{S\cup [\l]}}}}

\def\GK{{G_{K}}}

\def\linv{{1/\l}}

\def\I{{I_\l}}
\def\cI{{\mathcal{I}_\l}}

\def\Pminus{{\P^1 - \{0,1,\infty\}}}
\def\PminusC{{\P^1(\C) - \{0,1,\infty\}}}
\def\PminusQbar{{\P^1_\Qbar - \{0,1,\infty\}}}

\def\Hcts{H_{\mathrm{cts}}}
\def\Het{H_{\mathrm{\acute{e}t}}}

\def\Gm{{\mathbb{G}_m}}
\def\Ga{{\mathbb{G}_a}}
\def\Gml{{\Gm(\Ql)}}

\def\dot{{\bullet}}

\def\blank{\phantom{x}}

\def\MTMl{{{\mathcal M}TM_\l}}

\def\limproj#1{\lim_{\stackrel{\longleftarrow}{#1}}}
\def\liminj#1{\lim_{\stackrel{\longrightarrow}{#1}}}

\def\Homcts{\Hom^{\mathrm{cts}}}

\newcommand\im{\operatorname{im}}

\newcommand\Spec{\operatorname{Spec}}
\newcommand\Hom{\operatorname{Hom}}

\newcommand\Ext{\operatorname{Ext}}
\newcommand\Aut{\operatorname{Aut}}
\newcommand\Der{\operatorname{Der}}
\newcommand\Out{\operatorname{Out}}
\newcommand\Inn{\operatorname{Inn}}
\newcommand\OutDer{\operatorname{OutDer}}
\newcommand\Gal{\operatorname{Gal}}
\newcommand\Gr{\operatorname{Gr}}

\newcommand\cAut{\operatorname{{\mathcal Aut}}}
\newcommand\cOut{\operatorname{{\mathcal Out}}}

\newtheorem{theorem}{Theorem}[section]
\newtheorem{lemma}[theorem]{Lemma}
\newtheorem{proposition}[theorem]{Proposition}
\newtheorem{corollary}[theorem]{Corollary}

\theoremstyle{definition}
\newtheorem{definition}[theorem]{Definition}
\newtheorem{example}[theorem]{Example}
\newtheorem{conjecture}{Conjecture}

\theoremstyle{remark}
\newtheorem{remark}[theorem]{Remark}

\begin{document}

\title[Galois Actions on the Fundamental Group of $\Pminus$]
{Weighted Completion of Galois Groups and Galois Actions on
the Fundamental Group of $\Pminus$}

\author{Richard Hain}
\address{Department of Mathematics\\ Duke University\\
Durham, NC 27708-0320}
\email{hain@math.duke.edu}
\thanks{The first author was supported in part by grants from the National
Science Foundation. The second author was supported in part by a Mombusho
Grant and also by MSRI during a visit in the fall of 1999.}

\author{Makoto Matsumoto}
\address{Department of Mathematical Sciences\\ 
Faculty of Integrated Human Studies \\ Kyoto University\\
Yoshida-nihonmatsu-cho Sakyo-ku Kyoto 606-8501 Japan}
\email{matumoto@math.h.kyoto-u.ac.jp}

\date{\today}



\maketitle

\section{Introduction}

Fix a prime number $\l$. In this paper we prove a conjecture
\cite[p.~300]{Ihara3}, which Ihara attributes to Deligne, about
the action of the absolute Galois group on the pro-$\l$ completion of the
fundamental group of the thrice punctured projective line. It is stated
below.  Similar techniques are also used to prove part of a conjecture of
Goncharov \cite[Conj.~2.1]{goncharov}, also about the action of the absolute
Galois group on the fundamental group of the thrice punctured projective line,
and which derives from the conjectures of Deligne and Ihara and questions of
Drinfeld \cite[p.~859]{drinfeld}.

Ihara's version of Deligne's conjecture concerns the outer action
\begin{equation}
\label{outer_action}
\phi_\l : G_\Q \to \Out \pi_1(\PminusC,x)^\prol
\end{equation}
of the absolute Galois group $G_\Q$ on the pro-$\l$ fundamental group of
the thrice punctured projective line.
The pro-$\l$ completion of $\pi_1(\PminusC,x)$ is the inverse limit of its
finite quotients of $\l$-power order. We shall denote it by $\pi^\prol$.

Ihara defines a filtration
$$
G_\Q = \I^0 G_\Q \supseteq \I^1 G_\Q \supseteq \I^2 G_\Q \supseteq \cdots
$$
of $G_\Q$ by
\begin{equation}
\label{filtration}
\I^m G_\Q =
\ker\big\{\phi_\l : G_\Q \to \Out \big(\pi^\prol/L^{m+1}\pi^\prol\big)\big\}
\end{equation}
where $L^m \pi^\prol$ denotes the $m$th term of the lower central series
of $\pi^\prol$ indexed so that $L^1\pi^\prol = \pi^\prol$.
It is independent of the choice of base point $x$ but depends on the choice
of the prime $\l$, as the $\l$-adic cyclotomic character induces an
isomorphism
$$
\Gr_\I^0 G_\Q \cong \Zlx.
$$
It has the property that
$$
[\I^m G_\Q,\I^n G_\Q ] \subseteq \I^{m+n} G_\Q
$$
from which it follows that the positive part
$$
\Gr^{>0}_\I G_\Q := \bigoplus_{n>0} \Gr^n_\I G_\Q
$$
of the associated graded group of $G_\Q$ is a Lie algebra over $\Zl$ with
bracket induced by the group commutator. 

In this paper we shall prove:

\begin{conjecture}[{\cite[p.~300]{Ihara3}}]
\label{main}
The $\Ql$-form $(\Gr^{>0}_\I G_\Q)\otimes \Ql$
of the positive part of $\Gr^\dot_\I G_\Q$ is generated as a Lie algebra
by elements $s_3$, $s_5$, $s_7$, \dots where $s_{2n+1}\in\Gr^{2n+1}_\I G_\Q$.
\end{conjecture}

In fact, combining this result with the work of Ihara \cite{Ihara3}, one
can see (cf.\ Remark~\ref{h1_compn}) that there is a $G_\Q$-equivariant
isomorphism
$$
H_1(\Gr^{>0}_\I G_\Q,\Ql) \cong
\bigoplus_{n \ge 1} \Hcts^1(G_\Q,\Ql(2n+1))^\ast \otimes \Ql(2n+1)
\cong \bigoplus_{n\ge1}\Ql(2n+1).
$$
The value of Soul\'e's element \cite{soule} of $\Hcts^1(G_\Q,\Ql(2n+1))$
on the image of $s_{2n+1}$ in $\Hcts^1(G_\Q,\Ql(2n+1))^\ast\otimes
\Ql(2n+1)$ is non-zero.

In \cite{Ihara3} Ihara constructs explicit elements $\sigma_3, \sigma_5,\
\sigma_7, \dots$ of $\Gr^{>0}_\I G_\Q$, which are sometimes called
{\it Soul\'e elements}. He studies the non-vanishing of the Lie brackets of
these elements and asks \cite[p.~300]{Ihara3} whether the $\sigma_j$
generate $(\Gr^{>0}_\I G_\Q)\otimes \Ql$ freely. 
In \cite{Ihara9}, he proves that if $(\Gr^{>0}_\I G_\Q)\otimes \Ql$ is free,
then it is generated by the Soul\'e elements. He also shows that the
$\sigma_j$ generate an open subgroup of the image of the homomorphism
$$
\Gl \rightarrow \Out\big(\pi^\prol/L^{m+1}\pi^\prol\big)
$$
for each $m$, a consequence of Conjecture~\ref{main}. The techniques of this
paper do not seem to shed light on the question of freeness.

Our basic tool, which we develop in Paragraphs~\ref{neg_wtd_ext},
\ref{wt_filtn}, and \ref{wtd_cmp}, is the theory of weighted
completion of a profinite group $\G$ with respect to a continuous, Zariski
dense representation $\rho : \G \to R(\Ql)$ into a $\Ql$-group
endowed with a distinguished central cocharacter $w : \Gm \to R$. The
weighted completion in this case consists of
\begin{enumerate}
\item a proalgebraic $\Ql$-group $\cG$ which is an extension
$$
1 \to \U \to \cG \to R \to 1
$$
of $R$ by a prounipotent group $\U$ with the property that $H_1(\U)$,
when viewed as a $\Gm$-module via $w$, has only negative weights in the
sense of representation theory;
\item  a continuous homomorphism $\rhotilde : \G \to \cG(\Ql)$ which lifts
$\rho$.
\end{enumerate}
These data are required to be universal for homomorphisms of $\G$
into such ``negatively weighted extensions" of $R$ by a prounipotent group.
It is a variant of Deligne's notion of the relative Malcev completion, which
is developed in \cite{hain:comp,hain:malcev}. One can also view weighted
completions of $\G$ as the tannakian fundamental group of an appropriate
category of finite dimensional $\G$-modules 
over $\Ql$. This approach is explained
in \cite{hain-matsumoto:exp} --- see also Remark.~\ref{tannaka}.

Weighted completions arise naturally in Galois theory in many contexts, the
simplest of which is the following:
\begin{enumerate}
\item $\l$ is a prime number,
\item $\G = \Gl$, the Galois group of the maximal algebraic extension of $\Q$
unramified outside $\l$,
\item $R = \Gm$,
\item $\rho : \Gl \to \Gm(\Ql) = \Qlx$ is the $\l$-adic cyclotomic character,
\item $w : \Gm \to \Gm$ is the homomorphism defined by $w(x) = x^{-2}$.
\end{enumerate}
In Section~\ref{computn} we show that the weighted completion is an extension
$$
1 \to \Kl \to \Al \to \Gm \to 1
$$
where the Lie algebra $\kl$ of $\Kl$ is isomorphic (non-canonically) 
to the completion of the Lie algebra generated by
$$
\bigoplus_{n\ge 0} \Hcts^1(\Gl,\Ql(2n+1))^\ast\otimes \Ql(2n+1)
\cong \bigoplus_{n\ge 0}\Ql(2n+1).
$$

The most important property of weighted completions (or more generally,
{\it negatively weighted extensions}) is that each module of such a
group has a natural {\it weight filtration}, and that the weight graded
functor is exact on the category of such modules. This exactness property,
called {\it strictness}, is familiar from Deligne's mixed Hodge theory
\cite{deligne:hodge} as well as Galois theory. It is a key
ingredient in the proof of Conjecture~\ref{main}. The reason for taking the
central cocharacter to be $x \mapsto x^{-2}$ in the previous paragraph
is to make the representation theoretic weights agree with the usual
weights in algebraic geometry coming from Hodge theory and Galois theory.

The Lie algebra of a weighted completion is a module over the group
via the adjoint action and therefore has a natural weight filtration.
In particular, the Lie algebra $\al$ of the weighted completion of $\Gl$
has a weight filtration where $\kl = W_{-1}\al$. The copy of $\Ql(2n+1)$
in the previous paragraph lies in $W_{-2(2n+1)}\kl$ and projects to
a non-zero element in the $-2(2n+1)$th weight graded quotient of $\kl$.
These together generate $\Gr^W_\dot\kl$.

The first step in the proof of Conjecture~\ref{main} is to use the fact, due
to Ihara \cite{Ihara1}, that the Galois representation (\ref{outer_action})
factors through $\Gl$. This outer action induces one on the unipotent 
completion $\cP$ over $\Ql$ of the geometric fundamental group of $\Pminus$.
One can
define the filtration (\ref{filtration}) of $\Gl$ by considering this action.
The representation of $\Gl$ on the Lie algebra $\p$ of $\cP$ is negatively
weighted in an appropriate sense and induces a homomorphism
$\Al \to \Out \cP$ and a corresponding homomorphism $\al \to \OutDer \p$
on Lie algebras. Using strictness, one can show (see Section~\ref{actions})
that
$$
\big(\Gr_\I^m \Gl\big) \otimes_\Zl \Ql
\cong \Gr^W_{-2m}\im\{\kl \to \OutDer \p\}
\cong \im\big\{\Gr^W_{-2m} \kl \to \OutDer\Gr^W_\dot \p \big\}
$$
whenever $m>0$. Conjecture~\ref{main} then follows from the computation of
$\Gr^W_\dot\al$ in Section~\ref{computn}.

Goncharov's Conjecture is similar to Ihara's and concerns the action
(as distinct from the outer action) of the Galois group of the kernel
of the $\l$-adic cyclotomic character, $\Q(\mu_{\l^\infty})$, on the
$\l$-adic prounipotent completion $\cP$ of
$\pi_1(\PminusC,\overrightarrow{01})$
with base point the tangent vector at 0 pointing towards 1. In addition to
conjecturing that the Lie algebra of the Zariski closure of the image of
$G_{\Q(\mu_{\l^\infty})}$ in $\Aut \cP$ is generated by elements
$s_3, s_5, s_7, \dots$, he conjectures that it is freely generated by them.
We are able to prove the generation part of his conjecture, but not the
freeness. Goncharov's Conjecture is considered in Paragraph~\ref{goncharov}.

In Section~\ref{deligne} we prove $\l$-adic versions of
two related conjectures of Deligne, \cite[8.2, p.~163]{deligne} and
\cite[8.9.5, p.~168]{deligne}, concerning mixed Tate motives over the
spectrum of the ring of $S$-integers of a number field. As was
later pointed out to us by Goncharov, these had already been proved by
Beilinson and Deligne in their unpublished manuscript \cite{beilinson-deligne}.
Our approach via weighted completion is different, but equivalent, to the
approach of Beilinson and Deligne. Notably one should mention that Goncharov
\cite{goncharov:mtm} has used the work of Voevodsky \cite{voevodsky} and
Levine \cite{levine} to construct the category of mixed Tate motives over the
spectrum of a number field. It suggests a working definition
of the category of mixed Tate motives over the spectrum of a ring of
$S$-integers in a number field.

Our proof of Conjecture~\ref{main} can be conceptualized by saying that
the unipotent completion of the fundamental group of
$\pi_1(\Pminus,\overrightarrow{01})$ should be a mixed Tate motive over
$\Spec\Z$, and therefore a module over the tannakian fundamental group of
the category of mixed Tate motives over $\Spec \Z$. As was pointed out to
us by a referee, once one knows that the unipotent completion of the
fundamental group of $\pi_1(\Pminus,\overrightarrow{01})$ is a mixed Tate
motive over $\Spec \Z$ in the sense of Goncharov, Conjecture~\ref{main}
will follow by some arguments including the strictness of weight filtrations.
Our approach here is more direct, treating only $\l$-adic Galois 
representations.

Deligne's Conjectures concern
the category of mixed Tate motives over the spectrum $X_{F,S} :=
\Spec \O_F - S$ of the ring of $S$-integers of a number field $F$. 
His conjectures
essentially say that
\begin{enumerate}
\item there should be a category of mixed Tate motives over $X_{F,S}$;
\item the tannakian fundamental group of this category should be a 
proalgebraic $\Q$-group which is an extension of $\Gm$ by a prounipotent
group;
\item the Lie algebra of this prounipotent radical should be free and
its first cohomology is isomorphic to 
$\oplus_{n>0} K_n(X_{F,S})\otimes \Q(-n)$.
\end{enumerate}

Let $X_{F,S}[\linv]$ denote $X_{F,S}$ minus all primes over $\l$.
In Section~\ref{deligne}, we define the category of {\it $\l$-adic
mixed Tate modules} over $X_{F,S}$ to be the category whose objects are
$\pi_1(X_{F,S}[\linv])$-modules $M$ over $\Ql$ endowed with
a weight filtration $W_\dot$ by $\pi_1(X_{F,S}[\linv])$-submodules 
such that
all odd weight graded quotients of $M$ vanish and the $(-2m)$th graded
quotient is a sum of a finite number copies of $\Ql(m)$. In addition,
we require that all representations be crystalline at all primes $\p$
outside $S$ that lie over $\l$.

When $S$ contains all primes that lie over $\l$, the tannakian fundamental
group of the category of $\l$-adic mixed Tate modules is the weighted
completion $\AFS$
of $\pi_1(X_{F,S})$ with respect to the cyclotomic character. This is shown,
in Section~\ref{computn}, to be an extension of $\Gm$ by a prounipotent
$\Ql$-group whose Lie algebra is isomorphic to a completion of the free Lie
algebra generated by
$$
\bigoplus_{n>0} \Het^1(X_{F,S},\Ql(n))^\ast\otimes \Ql(n).
$$
This group is related to the algebraic $K$-theory of $X_{F,S}$ by the
regulators
$$
c_1 : K_{2n-1}(X_{F,S}) \to \Het^1(X_{F,S},\Ql(n)).
$$
Soul\'e's result \cite{soule} when $\l$ is odd, and results of 
Rognes and Weibel \cite{Rognes} when $\l = 2$, 
imply that these induce isomorphisms
$$
\bigoplus_{n>0} K_{2n-1}(X_{F,S})\otimes \Ql
\stackrel{\simeq}{\longrightarrow} \bigoplus_{n>0} \Het^1(X_{F,S},\Ql(n)).
$$
This yields a proof of an $\l$-adic version of Deligne's conjectures
similar to the unpublished proof of Beilinson and Deligne
\cite{beilinson-deligne}.

There is an appendix on unipotent completion where we collect and prove
results we need. Many of these must be well known but we know of no
good references, while other results may be new. One innovation is that
we introduce the notion of the $\l$-adic unipotent completion of a profinite
group and prove a comparison theorem which states that if $\G$ is a 
finitely generated group, then the $\l$-adic unipotent completion of 
the profinite (or pro-$\l$) completion of $\G$ is isomorphic to the
$\Ql$-form of the ordinary unipotent completion of $\G$.

There is a second appendix in which we prove some results about the 
continuous cohomology of the Galois groups $\GFS$ which are surely
well known to the experts, but for which we could not find suitable
references.
\bigskip

\noindent{\it Conventions.}
Unless mentioned otherwise, all fields in this paper will be of
characteristic zero.

The functor $H_1$ when applied to a group, Lie algebra, etc.\ will denote
the maximal abelian quotient. In the topological case, it will 
denote the quotient by the closure of the commutator subgroup or subalgebra,
as appropriate.

Following the conventions in Hodge theory (cf.\ \cite{deligne:hodge}), we
shall denote increasing filtrations with a subscript index:
$$
W_\dot M: \qquad \cdots \subseteq W_{m-1} M \subseteq W_m M \subseteq \cdots
$$
and decreasing filtrations with a superscript index:
$$
I^\dot M: \qquad \cdots \subseteq I^m M \subseteq I^{m+1} M \subseteq \cdots
$$
Their associated gradeds are defined by
$$
\Gr^W_\dot M := \bigoplus_{m\in \Z} W_mM/W_{m-1}M 
\text{ and } 
\Gr^\dot_I M := \bigoplus_{m\in \Z} I^m M/I^{m+1}M.
$$
Throughout this paper, $k$ denotes a field of characteristic zero.
\bigskip

\noindent{\it Acknowledgements.}
It is a pleasure to thank all those with whom we have had helpful
discussions: Akio Tamagawa, Masato Kurihara and Chad Schoen for discussions
on \'etale cohomology, Kazuya Kato and Takeshi Tsuji for discussions on
crystalline representations, and Yasutaka Ihara for his
generous comments on and clarification about the origin and history of the
conjectures considered in this paper.
We would also like to thank Alexandre Goncharov, for pointing out the
relevance of the work of Beilinson and Deligne \cite{beilinson-deligne},
Romyar Sharifi, for useful comments on and corrections to
the manuscript, and the referee who made many useful comments.

\section{Preliminaries on Proalgebraic Groups}
\label{sec:proalg}

In this paper, a group scheme or equivalently a $k$-group scheme
means an {\em affine} group scheme over $k$.
An algebraic group over a field $k$ or an algebraic $k$-group
means an {\em affine} group scheme of finite type over $k$.
A proalgebraic group is a projective limit of 
algebraic groups over $k$. It is known that the category of
proalgebraic groups over $k$ is equivalent to the category of 
affine group schemes over $k$, and any proalgebraic
group is a projective limit of a projective system 
consisting of surjective morphisms
(cf.\ \cite[Corollary 2.7, p.~128]{dmos}).

If $\displaystyle\cG=\limproj\alpha G_\alpha$, then the set of 
its $k$-rational points satisfies
$$
\cG(k) = \limproj \alpha G_\alpha(k).
$$
The Zariski closure of a subgroup $S$ of $\cG(K)$, 
$K$ an extension field of $k$,
is the smallest proalgebraic $k$-group that contains it --- it is the
projective limit $\varprojlim S_\alpha$ where $S_\alpha$ is the $k$-Zariski
closure of the image of $S$ in $G_\alpha(K)$.

The Lie algebra $\g$ of the proalgebraic
$k$-group
$$
\cG = \limproj \alpha G_\alpha,
$$
where each $G_\alpha$ is algebraic, is simply the inverse limit
$$
\limproj \alpha \g_\alpha
$$
of the Lie algebras of the $G_\alpha$. 
It will be regarded as a topological
Lie algebra with the topology defined by the inverse limit, where each
$\g_\alpha$ has the discrete topology.

If each $G_\alpha$ is unipotent, we say that $\cG$ is {\it prounipotent}.
The Lie algebra of a prounipotent group is the inverse limit of nilpotent
Lie algebras (in this paper a {\em nilpotent} Lie algebra means 
a {\em finite} dimensional nilpotent Lie algebra).  Such a Lie algebra is
said to be {\it pronilpotent}.

The well known correspondence between unipotent groups of $k$ and nilpotent
Lie algebras over $k$, which assigns to a unipotent group its Lie algebra,
extends easily to a correspondence between prounipotent
groups over $k$ and pronilpotent Lie algebras over $k$ (cf.\
\cite[Appendix~A]{quillen}).

A {\it graded Lie algebra} is a Lie algebra in the category of graded vector
spaces. It is of {\it finite type} if each of its graded quotients is finite
dimensional. The following results are needed in the sequel. The first is an
immediate consequence of the elementary fact that a homomorphism
$\n_1 \to \n_2$ between finite dimensional nilpotent Lie algebras is
surjective if and only if the induced map $H_1(\n_1) \to H_1(\n_2)$ is
surjective.

\begin{proposition}
\label{gen_lie}
If $\n = \oplus_{m>0} \n_m$ is a graded Lie algebra , then the image of any
section of the canonical surjection $\n \to H_1(\n)$ generates $\n$. \qed
\end{proposition}

Any such section induces a homomorphism from the free Lie algebra generated
by $H_1(\n)$. Thus we have:

\begin{corollary}
\label{free_quot}
If $\n = \oplus_{m>0} \n_m$ is a graded Lie algebra, then there is a free
graded Lie algebra $\f = \oplus_{m>0}\f_m$ and a surjective homomorphism
$\f \to \n$ such that the induced map $H_1(\f) \to H_1(\n)$ is an
isomorphism. \qed
\end{corollary}

Finally, we recall the Levi decomposition \cite[p.~158]{borel:GTM}.

\begin{proposition}
\label{split}
Every algebraic $k$-group $G$ is an extension of a reductive $k$-group
$R$ (called the Levi quotient) by a unipotent group $U$ (called the
unipotent radical),
$$
1 \to U \to G \to R \to 1.
$$
This extension is split and any two splittings differ by conjugation by
an element of $U(k)$. \qed
\end{proposition}

\section{Negatively Weighted Extensions and Weight Filtrations}
\label{neg_wtd_ext}

Let $R$ be an algebraic $k$-group.
Denote $\Gm_{/k}$ by $\Gm$. Suppose that $w : \Gm \to R$
is a central cocharacter --- that is, a homomorphism whose image is contained
in the center of $R$.

We shall denote by $Q^m$ the 1-dimensional irreducible representation of
$\Gm$ on which it acts by the $m$th power of the standard character.
We shall denote the $Q^m$-isotypical component of a $\Gm$-module $V$ by
$V_m$.

Suppose that
$$
1 \to U \to G \to R \to 1
$$
is an extension of $R$ by a unipotent group in the category of algebraic
$k$-groups. The first homology of 
$U$ is an $R$-module, and therefore a $\Gm$-module via the homomorphism
$w$. Thus we can write
$$
H_1(U) = \bigoplus_{m\in \Z} H_1(U)_m.
$$
We shall say that this extension is {\it negatively weighted} with 
respect to $w$ if $H_1(U)_m$ vanishes whenever $m\ge 0$.

We shall say that a proalgebraic group which is an extension of
$R$ by a prounipotent group $\U$ is {\it negatively weighted} if it is
an inverse limit of negatively weighted algebraic groups.

If $w$ is trivial, then there are no non-trivial negatively weighted
extensions of $R$.

\subsection{The Weight Filtration}
\label{wt_filtn}
Suppose that an algebraic $k$-group $G$ is a negatively weighted extension of
$R$ with respect to the central cocharacter $w : \Gm \to R$.
In this section, we show that representations
of a negatively weighted extension of $R$ have a natural weight filtration,
and that morphisms between such modules are strict with respect to this weight
filtration.

\begin{lemma}
\label{lifts}
There is a lift of the homomorphism $w : \Gm \to R$ to a homomorphism
$\wtilde : \Gm \to G$, and any two such liftings are conjugate by an element
of $U(k)$.
\end{lemma}

\begin{proof}
This can be seen by pulling back the extension $U \to G \to R$
along $w$ and then applying Proposition~\ref{split}.
\end{proof}

Fix such a splitting $\wtilde$. Each finite dimensional $G$-module $V$ can
then be regarded as a $\Gm$-module via $\wtilde$, and can therefore be
decomposed
\begin{equation}
\label{weight_split}
V = \bigoplus_{n\in \Z} V_n
\end{equation}
under the $\Gm$-action. Define the {\it weight filtration} of $V$ by
$$
W_n V = \bigoplus_{m \le n} V_m.
$$
A priori, this filtration depends on $\wtilde$, but we will show, in 
Proposition~\ref{weight_indep}, that it does not.

The $n$th weight graded quotient of $V$ is defined by
$$
\Gr^W_n V := W_n V / W_{n-1} V.
$$
The inclusion $V_n \hookrightarrow V$ induces a natural isomorphism
$$
\Gr^W_n V \cong V_n.
$$

It follows directly from the representation theory of $\Gm$ that, for a
fixed splitting $\wtilde$, the weight splitting (\ref{weight_split}) is 
compatible with Hom and tensor products.

\begin{proposition}
If $U$ and $V$ are $G$-modules, then
$$
\Hom(U,V)_n = \bigoplus_{r-m=n} \Hom(U_m,V_r)
$$ 
and
$$
(U\otimes V)_n = \bigoplus_{r+m=n} U_m\otimes V_r.
$$
Moreover, if $\phi :  V \to W$ is $G$-equivariant, then $f(V_n) \subset W_n$
for all $n$. \qed
\end{proposition}

An immediate consequence is that the weight filtration is compatible with
respect to Hom and tensor products.

\begin{corollary}
\label{cor:homfiltration}
If $U$ and $V$ are $G$-modules, then
$$
W_n \Hom(U,V) = \{\phi \in \Hom(U,V) : \phi(W_r U) \subseteq W_{r+n} V
\mbox{ for all $r\in \Z$}
\}
$$
and
$$
W_n(U\otimes V) = \sum_{r+s = n} W_rU \otimes W_sV. \qed
$$
\end{corollary}

We can apply this to the adjoint action to decompose $\g$, the Lie algebra
of $G$, and $\u$, the Lie algebra of the negative unipotent part $U$:
$$
\g = \bigoplus_{n\in \Z} \g_n \text{ and } \u = \bigoplus_{n\in \Z} \u_n.
$$
Denote the Lie algebra of $R$ by $\r$.

\begin{proposition} With notation as above,
\begin{enumerate}
\item if $V$ is a $G$-module, $x\in \g_n$ and $v \in V_m$, then
$x\cdot v \in V_{n+m}$;
\item $\g$ is a graded Lie algebra; that is, if $x \in \g_n$ and $y\in \g_m$,
then $[x,y] \in \g_{n+m}$;
\item $\u$ is a sub graded Lie algebra of $\g$, that is, $\u_m \subseteq \g_m$;
\item if $n \ge 0$, then $\u_n= 0$;
\item the graded quotients of $\g$ are given by
$$
\g_n = 
\begin{cases}
0 & n > 0; \cr
\r & n = 0; \cr
\u_n & n<0.
\end{cases}
$$
\end{enumerate}
\end{proposition}

\begin{proof}
The first statement follows from the previous proposition as the mapping
$\g \otimes V \to V$ is $G$-equivariant. The second follows from the
first by taking $V$ to be the adjoint representation. The third statement
is clear.

Set
$$
\u_{<0} = \bigoplus_{n<0} \u_n.
$$
The inclusion of this into $\u$ induces a surjection on $H_1$ as 
$H_1(\u)$ has all negative weights.  Since $\u$ and $\u_{<0}$ are both
nilpotent, this implies that they are equal.  This proves the fourth
statement.

The exact sequence
$$
0 \to \u \to \g \to \r \to 0
$$
is an exact sequence of graded $\g$-modules. Since $w$ is central in $R$,
it follows that $\r = \r_0$. Since $\u$ has only negative weights, it
follows that the projection $\g \to \r$ induces an isomorphism
$\g_0 \cong \r$ and that $\g_n = \u_n$ whenever $n \neq 0$.
\end{proof}

\begin{corollary}
With notation as above, $W_\dot \g$ is a filtration of $\g$ by Lie 
ideals, and we have
$\g = W_0 \g$, $\u = W_{-1} \g$, and $\r \cong \Gr^W_0 \g$. \qed
\end{corollary}

\begin{lemma}
\label{u_pres}
For every choice of lift $\wtilde$, each term $W_n V$ of the weight
filtration of the $G$-module $V$ is a $\g$-module, and $U$ acts
trivially on its associated graded $\Gr^W_\dot V$.
\end{lemma}

\begin{proof}
Since $\u = W_{-1} \g$ and since the action $\g \otimes V \to V$ is
compatible with the gradings, each term $W_n V$ of the weight filtration
is a $\g$-module. Since $\u = W_{-1}\g$, the image of the action
$\u \otimes W_n V \to W_n V$ is contained in $W_{n-1} V$, from which the
result follows.
\end{proof}

The Levi decomposition extends to negatively weighted extensions.

\begin{lemma}
\label{wtd_levi}
If $1 \to U \to G \to R \to 1$ is a negatively weighted extension with
respect to the central cocharacter $w:\Gm \to R$, then the projection
$G \to R$ splits, and any two splittings are conjugate by an element of
$U$.
\end{lemma}

\begin{proof}
If $R$ is reductive, this follows directly from the classical Levi 
decomposition (Prop.~\ref{split}). When $R$ connected, the existence
of the splitting follows from the existence of a splitting of $\g \to \r$
as $U$ is simply connected. The existence of such a Lie algebra splitting
follows by first taking a lifting $\wtilde : \Gm \to G$ of $w$, and then
decomposing $\g$ under the restriction of the adjoint action to $\Gm$ via
$\wtilde$. The summand of weight zero (the invariant part) is a Lie algebra
(as the bracket preserves weights) and projects isomorphically to $\r$.

First we will prove that any two liftings have to be conjugate by an
element of $U$. The intersection of $U$ with the centralizer $C(\wtilde)$
of a lift $\wtilde : \Gm \to G$ of $w$ is just the identity as the Lie
algebra of $U$ has strictly negative weights with respect to $\wtilde$ under
the adjoint action. It follows that if $s : R \to G$ is a splitting and
$\wtilde = s\circ w$, then $s(R) = C(\wtilde)$. Since the lifts $\wtilde$
are unique up to conjugation by an element of $U$, any two splitting of
$G \to R$ have to be conjugate by an element of $U$.

To prove the existence of the splitting in general, first take a Levi
decomposition
$$
R \cong \Rbar \ltimes W
$$
of $R$, where $\Rbar$ is reductive and $W$ unipotent. This also induces a
Levi decomposition
$$
R^o \cong \Rbar^o \ltimes W
$$
of $R^o$, the connected component of the identity. By the connected case,
we can lift $R^o$ to $G$. In particular, we can lift $W$ to $G$. Since
$W$ and $U$ are unipotent, $\Rbar$ is the reductive quotient of $G$. By
the classical Levi decomposition, we have a splitting of $G \to \Rbar$.
The subtlety is that we have to choose this splitting so that its image
normalizes the lift of $W$.

The splittings $\Rbar \to G$ and $R^o \to G$ determine lifts $\wtilde_0$
and $\wtilde_1$ of $w$. By conjugating one of the splittings by an element
of $U$, we may assume they are the same. But then the images of both
splittings lie in the centralizer $C(\wtilde_0)$. It follows that the
image projection $C(\wtilde_0) \to R$ is surjective, and therefore an
isomorphism. Inverting this isomorphism provides the splitting.
\end{proof}

\begin{proposition}\label{weight_indep}
The weight filtration of a $G$-module $V$ is independent of the choice of
the lift $\wtilde$ of $w$ and has the property that each $W_m V$ is a
$G$-module. Moreover, each $\Gr^W_m V$ is an $R$-module of weight $m$.
\end{proposition}

\begin{proof}
Suppose that $\wtilde':\Gm \to G$ is another lift of $w$. Then, by
Lemma~\ref{lifts},
there is $u \in U(k)$ such that $\wtilde' = u\wtilde u^{-1}$. Write
$$
V = \bigoplus_{n\in \Z} V_n'
$$
where $V_n'$ is the subspace of $V$ where $\Gm$ acts with weight $n$
via $\wtilde'$. Then $V_n' = uV_n$. Since each term of the weight
filtrations associated to $w$ and $\wtilde$ are $G$-modules, it follows that,
for each $n$,
$$
\bigoplus_{m\le n} V_m' = \bigoplus_{m\le n} u V_m
= u \left(\bigoplus_{m\le n} V_m\right) \subseteq \bigoplus_{m\le n} V_m.
$$
The reverse inclusion follows by reversing the roles of $w$ and $\wtilde$.
The equality of the weight filtrations follows.

To prove the second assertion, we choose a splitting $s : R \to G$ of the
canonical projection $G \to R$, whose existence is assured by
Lemma~\ref{wtd_levi}. Choose the lift $\wtilde$ to be the composite
of $w: \Gm \to R$ with the splitting $s$. Let $V = \oplus V_n$ be the
corresponding weight decomposition of $V$. Since $w$ is central in $R$, it
follows that the image of $s$ preserves this decomposition and that each
$W_m V$ is an $(\im s)$-module. The result follows from Lemma~\ref{u_pres}.

The last assertion is clear.
\end{proof}

\subsection{The Proalgebraic Case}
Suppose now that $\cG$ is a proalgebraic $k$-group which is a negatively
weighted extension of the algebraic $k$-group $R$ and that $V$ is a finite
dimensional $\cG$-module. Note that $\cG$ is the projective limit of algebraic 
$k$-groups $G_\alpha$, each of which is a (necessarily negatively weighted)
extension of $R$. Since $V$ is finite dimensional, it is a $G_\alpha$-module
for some $\alpha$, and therefore has a weight filtration. It is easily
seen that this weight filtration is independent of $\alpha$. We therefore
have:

\begin{theorem}
\label{xx}
If $\cG$ is a proalgebraic $k$-group which is a negatively weighted
extension of the algebraic $k$-group $R$, then every finite dimensional
$\cG$-module has a natural weight filtration $W_\dot$ by $\cG$-submodules.
This has the following properties:
\begin{enumerate}
\item the action of $\cG$ on $\Gr^W_\dot V$ factors through $\cG \to R$;
\item the action of $\Gm$ on $\Gr^W_m V$ via $w$ is of weight $m$;
\item the weight filtration is compatible with $\Hom$ and $\otimes$
in that it obeys the rules in Corollary~\ref{cor:homfiltration}. \qed
\end{enumerate}
\end{theorem}

The weight filtration extends naturally to two different classes of
infinite dimensional representations of $\cG$ --- viz, those representations
that are projective or inductive limits of finite dimensional $\cG$-modules.
In particular, the Lie algebras $\g$ of $\cG$ and $\u$ of $\U$ have natural
weight filtrations.

\subsection{Morphisms}
\label{homom}
The weight filtration has nice naturality and exactness properties. The
following setting is a little technical, but convenient for our applications.

Suppose that $R_1$ and $R_2$ are algebraic $k$-groups with distinguished
central cocharacters $w_i : \Gm \to R_i$ and that $\cG_i$ is a proalgebraic
$k$-group where $i=1,2$. Suppose that $\cG_i \to R_i$ is an extension (not
necessarily negatively weighted) of $R_i$, where $i=1,2$. Suppose that
$\Phi : \cG_1 \to \cG_2$ and $\phi: R_1 \to R_2$ are homomorphisms
such that the diagram
$$
\begin{CD}
\cG_1 @>>> R_1 @<{w_1}<< \Gm \cr
@V{\Phi}VV @V{\phi}VV @| \cr
\cG_2 @>>> R_2 @<{w_2}<< \Gm \cr
\end{CD}
$$
commutes.

The following result says that, in some sense, negatively weighted
extensions are closed under subs and quotients. It is easily proved.

\begin{lemma}
\label{lem:subgp}
Suppose that $\Phi : \cG_1 \to \cG_2$ is as above.
\begin{enumerate}
\item If $\Phi$ is an inclusion and $\cG_2$ is a negatively weighted 
extension of $R_2$ with respect to $w_2$, then $\cG_1$ is a negatively
weighted extension of $R_1$ with respect to $w_1$;
\item If $\Phi$ is a surjection and $\cG_1$ is a negatively weighted 
extension of $R_1$ with respect to $w_1$, then $\cG_2$ is a negatively
weighted extension of $R_2$ with respect to $w_2$. \qed
\end{enumerate}
\end{lemma}

A $\Phi$-homomorphism from a $\cG_1$-module $V_1$ to a $\cG_2$-module $V_2$
is a $\cG_1$ homomorphism, where $V_2$ is regarded as a $\cG_1$-module via
$\Phi$.

We can form a category: The objects are 4-tuples $(V,R,w,\cG)$ where $R$ is an
algebraic $k$-group, $w:\Gm \to R$ a central cocharacter, $\cG$ a proalgebraic
group which is a negative extension of $R$, and where $V$ is a finite
dimensional $\cG$-module. Morphisms between such modules are the $\Phi$
homomorphisms defined above. We shall call this the category of {\it weighted
modules}. Note that every object of this category has a natural weight
filtration.

\subsection{Strictness}
Strictness is a key ingredient in our proof of Conjecture~\ref{main}.
A linear mapping $f: (V_1,W_\dot) \to (V_2,W_\dot)$ between two filtered
vector spaces is said to be {\it strict} with respect to the filtrations
$W_\dot$ if it is filtration preserving and if
$$
\im f \cap W_mV_2 = f(W_m V_1)
$$
for all $m\in \Z$.

A filtration $W_\dot$ of a vector space induces a filtration on each
subspace $A \hookrightarrow V$ and each quotient $q:V \twoheadrightarrow B$
by
$$
W_mA := A \cap W_m V \text{ and } W_mB = q(W_m V).
$$

\begin{proposition}
If $f: (V_1,W_\dot) \to (V_2,W_\dot)$ is strict with respect to $W_\dot$,
then there are natural isomorphisms
$$
\ker \Gr^W_\dot f \cong \Gr^W_\dot \ker f \text{ and }
\im \Gr^W_\dot f \cong \Gr^W_\dot \im f
$$
of graded vector spaces. \qed
\end{proposition}

\begin{theorem}
\label{naturality}
Each homomorphism $f : V_1 \to V_2$ in the category of weighted modules
preserves the weight filtration and is strict with respect to it.
Consequently, the functors
$$
V \mapsto \Gr^W_m V, \quad V \mapsto W_m V, \text{ and } V \mapsto V/W_mV
$$
are all exact on the category of finite dimensional $\cG$-modules.
\end{theorem}

\begin{proof}
We use the notation of Paragraph~\ref{homom}. Since $V_1$ and $V_2$
are finite dimensional, we may assume that $\cG_1$ and $\cG_2$ are
algebraic. To prove that the weight
filtration is preserved by the morphism $f:V_1 \to V_2$, choose a
lift $\wtilde_1 : \Gm \to G_1$, and define $\wtilde_2$ to be
$\Phi\circ \wtilde_1$. Then both $V_1$ and $V_2$ become $\Gm$-modules,
and $f$ is $\Gm$-equivariant. The first assertion follows. Strictness
follows similarly as $f$ becomes a map of graded vector spaces after
choosing the lift $\wtilde$.

The exactness of $\Gr^W_n$ follows for similar reasons. If
$$
0 \to V' \to V \to V'' \to 0
$$
is an exact sequence of weighted modules, one can choose compatible liftings
of the central cocharacters, so that the exact sequence becomes an exact
sequence of $\Gm$-modules. The sequence
$$
0 \to V_n' \to V_n \to V_n'' \to 0
$$
is thus exact for all $n$. The exactness of $\Gr^W_n$ follows as the $n$th
weight graded quotient of a module is naturally isomorphic to its weight $n$
part.
\end{proof}

This theorem extends to certain infinite dimensional
representations of such groups. There are two appropriate categories of such
$\cG$-modules. Namely, the category of projective limits of finite
dimensional representations of $\cG$, and the category of inductive limits of
finite dimensional representations of $\cG$. In both cases, there is a weight
filtration and morphisms are strict with respect to it.

\section{Weighted Completion}
\label{wtd_cmp}

We will work only in the profinite setting as that is where our examples
lie, although there is an obvious analogue in the discrete case. Denote
$\Gm_{/\Ql}$ by $\Gm$. Suppose that:
\begin{enumerate}
\item $\G$ is a profinite group;
\item $R$ is an algebraic group defined over $\Ql$;
\item $w : \Gm \to R$ is a central cocharacter;
\item $\rho : \G \to R(\Ql)$ is a continuous homomorphism with Zariski dense
image.
\end{enumerate}

\begin{definition}\label{def:weightedcompletion}
The {\it weighted completion of $\G$ with respect to $\rho$ and $w$}
is a pro-algebraic $\Ql$-group $\cG$, which is an extension of $R$ by a
prounipotent group $\U$, and a continuous homomorphism
$\rhotilde: \G \to \cG(\Ql)$ which lifts $\rho$:
$$
\begin{CD}
1 @>>> \U(\Ql) @ >>> \cG(\Ql) @>>> R(\Ql) @>>> 1\cr
@. @. @A{\rhotilde}AA @AA{\rho}A \cr
@. @. \G @= \G
\end{CD}
$$
It is characterized by its universal mapping property: if $G$ is a
proalgebraic $\Ql$-group which is a {\em negatively weighted} extension of
$R$ (with respect to $w$)
by a prounipotent group, and if $\phi : \G \to G(\Ql)$ is a homomorphism that
lifts $\rho$, then there is a unique homomorphism of proalgebraic groups
$\Phi : \cG \to G$ that commutes with the projections to $R$ and such that
$\phi = \Phi\rhotilde$.
\end{definition}

\begin{definition}\label{def:wtG-mod}
A $\Ql\G$-module $V$ is a $\Ql$-vector space with a continuous $\G$-action. 
A {\it weighted $\G$-module} is a finite dimensional $\Ql\G$-modules $V$
with an increasing filtration $W_\dot$ by $\G$-submodules where, for each
$n$, the action of $\G$ on $\Gr^W_n V$ factors through an action of $R$ via
$\rho$ and where
$$
\Gm \stackrel{w}{\to} R \to \Aut \Gr^W_n V
$$
is the $n$th power of the standard character.

Morphisms between weighted $\G$-modules are those as $\G$-modules.
\end{definition}

\begin{remark}
\label{tannaka}
The weighted completion of $\G$ with respect to $\rho : \G \to R(\Ql)$
can also be defined to be the tannakian fundamental group of the category
of weighted $\G$-modules. This is proved in \cite{hain-matsumoto:exp}.
\end{remark}

\begin{proposition}\label{prop:catconstruction}
The weighted completion of $\G$ with respect to $\rho : \G \to R(\Ql)$
and $w : \Gm \to R$ always exists.
\end{proposition}

\begin{proof}
Consider the category whose objects are pairs
$$
(G \to R,\phi : \G \to G(\Ql)), 
$$
where $G$ is an algebraic $\Ql$-group which is a negatively weighted extension of $R$
and where $\phi$ is a continuous, Zariski dense representation which
lifts $\rho : \G \to R(\Ql)$. A morphism
$$
F: (G_1 \to R,\phi_1 : \G \to G_1(\Ql))
\to (G_2 \to R,\phi_2 : \G \to G_2(\Ql))
$$
consists of a homomorphism $f : G_1 \to G_2$ of $\Ql$-groups that is compatible
with the projections to $R$ and satisfies $\phi_2 = f\cdot \phi_1$.

Since the image of $\phi$ is Zariski dense, there is at most one
morphism between any two objects. In addition, there is fibered product:
$$
(G_1 \to R,\phi_1 : \G \to G_1(\Ql)) \times_R
(G_2 \to R,\phi_2 : \G \to G_2(\Ql))
$$
is defined to be
$$
(G \to R,\phi : \G \to G(\Ql))
$$ 
where $G$ is the Zariski closure of the image of
$$
\phi_1\times_R \phi_2 : \G \to (G_1 \times_R G_2)(\Ql).
$$
It follows that this category is a projective system. The weighted completion
of $\G$ with respect to $\rho$ is then $(\cG,\rhohat : \G \to \cG(\Ql))$,
the projective limit of all objects of this category.
\end{proof}

\subsection{Cohomology of $\g$ and $\u$}
We shall view $\g$ and $\u$ as topological Lie algebras where neighbourhood
of zero are the kernels $\n$ of the canonical homomorphisms to the finite
dimensional quotients of $\g$ that are extensions
$$
0 \to \u/\n \to \g/\n \to \r \to 0.
$$
of $\r$ by a nilpotent Lie algebra.

Suppose that $V$ is a finite dimensional $R$-module. Define the continuous
cohomology of $\u$ and $\g$ with coefficients in $V$ to be the direct limit
of the cohomology of their canonical finite dimensional quotients:
$$
\Hcts^\dot(\u) = \liminj \n H^\dot(\u/\n)
$$
and
$$
\Hcts^\dot(\g,V) = \liminj \n H^\dot(\g/\n,V)\text{ and }
\Hcts^\dot(\u,V) = \liminj \n H^\dot(\u/\n,V) = \Hcts^\dot(\u)\otimes V.
$$
These can be computed using continuous cochains; the group $\Hcts^\dot(\g,V)$
is the cohomology of the complex
$$
\Homcts(\Lambda^\dot \g ,V) := \liminj \n \Hom(\Lambda^\dot(\g/\n),V)
$$
of continuous cochains. The differential is induced by the dual of the
bracket. A complex which computes $\Hcts^\dot(\u,V)$ is obtained by
replacing $\g$ by $\u$.

Since the Lie algebras $\g$ of $\cG$ and $\u$ of $\U$ are inverse limits
of the Lie algebras of the finite dimensional quotients of $\cG$, it follows
from Theorem~\ref{naturality} that each has a natural weight filtration. 
Since the category of $\cG$-modules is abelian, these pass to
$\Hcts^\dot(\u)$.

\begin{proposition}
\label{weights}
The Lie algebras $\g$ of $\cG$, $\r$ of $R$ and $\u$ of $\U$ have natural
weight filtrations. These satisfy
$$
\g = W_0 \g,\quad W_n\u = W_n\g \text{ whenever } n<0.
$$
These weight filtrations pass to cohomology. In particular, each cohomology
group $\Hcts^m(\u)$ has a canonical weight filtration which satisfies
$$
\Gr^W_m \Hcts^n(\u) = 0 \text{ if } m < n.
$$
More precisely, if $W_{m-1} \Hcts^1(\u) = 0$, then $W_{nm-1}\Hcts^n(\u) = 0$.
\end{proposition}

\begin{proof}
We need only prove the last statement. Since $\Hcts^1(\u)$ has weights
$\ge m$, $\u$ has weights $\le -m$, and its continuous dual $\u^\ast$
(a direct limit of finite dimensional $\cG$-modules) has weights $\ge m$.
It follows that the space
$$
\Homcts(\Lambda^\dot \u,k)
$$
of degree $n$ continuous cochains on $\u$ has weights $\ge nm$. Since the
bracket is a morphism, the Chevalley-Eilenberg complex
$$
(\Homcts(\Lambda^\dot \u,k), -[\blank,\blank]^\ast)
$$
of continuous cochains on $\u$ is a complex of $\cG$-modules. Its cohomology
is therefore a graded $\cG$-module and thus has a weight filtration induced
from that of $\u$. Since the  $n$-cochains have weights $\ge nm$, it follows
that the weights on $H^n(\u)$ are also $\ge nm$.
\end{proof}

\subsection{Basic structure of $\g$ and $\u$}

Each finite dimensional representation $V$ of $R$ can be decomposed
$$
V = \bigoplus_{n\in \Z} V_n
$$
under the $\Gm$ action. We shall call $V$ {\it pure} of weight $n$ if
$V=V_n$, and {\it negatively weighted} if $V_n = 0$ for all $n \ge 0$.
Since $w$ is central, each $V_n$ is an $R$-module.
Thus, if $V$ is indecomposable (i.e., cannot be written as a direct sum of
two non-trivial submodules), $V$ is pure of some weight we shall denote by
$n(V)$.

Denote by $\Hcts^i(\G, V)$ the cohomology of the complex of {\em continuous}
cochains
$$
f : \G^{m+1} \to V
$$
where $V$ is viewed as a continuous $\G$-module via $\rho : \G \to R(\Ql)$.
This is the same notion of continuous cohomology as used in
\cite[Sect.~2]{tate}.

\begin{theorem}
\label{h1_new}
For all pure $R$-modules $V$, there are natural isomorphisms
$$
\Homcts_R(H_1(\U),V) \cong
\begin{cases}
\Hcts^1(\G,V) & \text{ when } n(V) < 0; \cr
0 & \text{ when } n(V) \ge 0
\end{cases}
$$
and
$$
\Homcts_R(\Gr^W_n H_1(\U),V) \cong
\begin{cases}
\Hcts^1(\G,V) & \text{ when } n = n(V) < 0; \cr
0 & \text{ otherwise.}
\end{cases}
$$
\end{theorem}

The second assertion follows directly from the first, which is
proved in Section~\ref{proof:h1_new}.

\begin{corollary}\label{ident}
If $\Hcts^1(\G,V)$ vanishes for all pure representations
$V$ of $R$ of negative weight, then the weighted completion of
$\G$ relative to $\rho : \G \to R$ is $\rho$ itself.
\end{corollary}

\begin{proof}
The theorem above implies that $H_1(\u)=0$. Since $\u$ is pronilpotent,
$\u$ is trivial.
\end{proof}

\subsection{The case $R$ reductive}
\label{reductive}
Suppose now that $R$ is reductive. Let $\{V_\alpha\}_\alpha$ be a set of
representatives of the isomorphism classes of finite dimensional
irreducible representations of $R$. For convenience, we set
$n(\alpha) = n(V_\alpha)$.

The next result follows directly from Theorem~\ref{h1_new}.

\begin{theorem}
\label{h1}
If $\Hcts^1(\Gamma, V_\alpha)$ is finite dimensional for all $\alpha$ with
$n(\alpha)<0$, then
$$
\Hcts^1(\u)\cong \bigoplus_{\{\alpha:n(\alpha)\ge 1\}}
\Hcts^1(\G,V_\alpha^\ast)\otimes V_\alpha
$$
and
$$
H_1(\U) \cong H_1(\u) \cong
\limproj F
\bigoplus_{\alpha \in F} \Hcts^1(\G,V_\alpha)^\ast \otimes V_\alpha
$$
where $F$ ranges over the finite subsets of $\{\alpha:n(\alpha)<0\}$.
More precisely, $W_0 \Hcts^1(\u)$ vanishes and, if $m > 0$, then
$$
\Gr^W_m \Hcts^1(\u) \cong \bigoplus_{\{\alpha:n(\alpha) = m\}}
\Hcts^1(\G,V_\alpha^\ast)\otimes V_\alpha. \qed
$$
\end{theorem}

In the reductive case, we are able to control $\Hcts^2(\u)$ as well, which
will allow us to control the relations in $\u$.

\begin{theorem}
\label{h2}
There is a natural injective $R$-invariant homomorphism
$$
\Phi : \Hcts^2(\u) \hookrightarrow
\bigoplus_{\{\alpha:n(\alpha)\ge 2\}}
\Hcts^2(\G,V_\alpha^\ast)\otimes V_\alpha.
$$
In addition, if $m > 1$, then
$$
\Gr^W_m \Hcts^2(\u) \subseteq \bigoplus_{\{\alpha:n(\alpha) = m\}}
\Hcts^2(\G,V_\alpha^\ast)\otimes V_\alpha.
$$
\end{theorem}

The proof is given in Section~\ref{proof:h2}.

\begin{corollary}
\label{freeness}
Suppose that $\Hcts^1(\G,V_\alpha^\ast)$ vanishes for all $\alpha$ satisfying
$0 < n(\alpha) < d$. If $\Hcts^2(\G,V_\alpha^\ast)$ vanishes for all $\alpha$
with $n(\alpha) \geq 2d$, then $\Gr^W_\dot\u$ is a free Lie algebra.
Moreover, any lift of a basis of $\Gr^W_\dot H_1(\u)$ to $\Gr^W_\dot\u$
is a free generating set.
\end{corollary}

\begin{proof}
The condition on $\Hcts^1(\Gamma,V_\alpha^\ast)$ implies, by Theorem~\ref{h2},
that $W_{d-1}\Hcts^1(\u)=0$.
The last assertion in Proposition~\ref{weights} implies that
$W_{2d-1}\Hcts^2(\u) = 0$. The assumption about $\Hcts^2(\Gamma,V_\alpha^\ast)$
and the last assertion of Theorem~\ref{h2} imply that
$$
\Gr^W_m \Hcts^2(\u) = 0
$$
when $m\ge 2d$, so that $\Gr^W_\dot\Hcts^2(\u) = 0$. Since $\Gr^W_\dot$ is
an exact functor, it commutes with homology. So it follows that
$H^2(\Gr^W_\dot\u)$ vanishes. The result now follows from the following lemma.
\end{proof}

\begin{lemma}
\label{free_crit}
If $\n = \oplus_{m<0} \n_m$ is a graded Lie algebra over a field of
characteristic zero, then $\n$ is free if and only if $H^2(\n) = 0$.
\end{lemma}

\begin{proof}
By Corollary~\ref{free_quot}, there is a free graded Lie algebra
$\f = \oplus_{m<0} \f_m$ and a graded Lie algebra homomorphism
$\f \to \n$ which is surjective and induces an isomorphism on $H_1$.
Denote the kernel of this by $\a$. Note that $\n$ is free if and only
if $\a = 0$. Since $\a$ is an ideal in a negatively graded Lie algebra,
$\a = 0$ if and only if $\a/[\f,\a] = 0$.

There is a spectral sequence
$$
E_2^{a,b} = H^a(\n,H^b(\a)) \Rightarrow H^{a+b}(\f).
$$
Since $\f$ is free, $H^2(\f) = 0$. It follows that
$$
d_2 : H^0(\n,H^1(\a)) \to H^2(\n)
$$
is an isomorphism. The result follows as $H^0(\n,H^1(\a))$ is the graded
dual of $\a/[\f,\a]$.
\end{proof}

\begin{example}
\label{zlx}
Suppose that $\G = \Zlx$, that $R=\Gm_{/\Ql}$ and that
$\rho : \Zlx \hookrightarrow \Gm(\Ql)=\Qlx$ is the natural inclusion.
Take $w$ to be the inverse of the square of the standard character.%
\footnote{We use the inverse of 
the square of the standard character so that later, when
considering the weighted completion of $\pi_1(\Spec\Z[1/\l])$, for example,
the representation theoretic weights agree with weights defined using 
Frobenius. This square makes no difference to the weighted completion,
it just {\em doubles} the indexing of the weight filtration.}
In this example we compute the weighted completion of $\Zlx$ with respect
to $\rho$ and $w$. Note that
$$
\Hcts^1(\Zlx,Q^m) = 0,
$$
for all non-zero $m\in \Z$. This can be seen as follows.

For each $a \in \Zlx$, the mapping
$$
a : Q^m \to Q^m \quad x \mapsto a^m x
$$
is $\Zlx$-equivariant, and thus induces an automorphism of $H^n(\Zlx,Q^m)$.
A general fact in group cohomology implies that this automorphism is trivial.
(cf.\ \cite[p.~116]{serre:localfield}.) It follows that $(1-a^m)$ annihilates
$H^n(\Zlx,Q^m)$ for all $a\in \Zlx$. But if $m\neq 0$ and $a\neq 1$, then
$(1-a^m)$ is in $\Qlx$. It follows that $H^n(\Zlx,Q^m)$ vanishes when
$m\neq 0$.

Theorem~\ref{h1} tells us
that the unipotent radical $\U$ of the weighted completion of $\Zlx$ is
trivial, so that the weighted completion of $\Zlx$ with respect
to $\rho$ is just $\rho : \Zlx \to \Gm(\Ql)$.
\end{example}

More generally, we can take $\G$ in the previous example to be any open
subgroup of $\Zlx$. One has the same vanishing of cohomology and
consequently, that the weighted completion of $\G$ with respect to
$\rho : \G \to \Gm(\Ql)$ is just $\rho$ itself.

\subsection{Naturality}
\label{nat:setup}
Here we record some more technical naturality statements that are needed
in the sequel.

Suppose that $f : \G_1\to \G_2$ is a continuous homomorphism between two
profinite groups. Suppose that $R_1$ and $R_2$ are two algebraic $\Ql$-%
groups with central cocharacters $w_1$ and $w_2$, respectively. Suppose that
$$
\rho_1 : \G_1 \to R_1(\Ql) \text{ and } \rho_2 : \G_2 \to R_2(\Ql)
$$
are two continuous, Zariski dense homomorphisms. We can form the
weighted completions $\cG_1$ and $\cG_2$ of $\G_1$ and $\G_2$. Suppose that
$\phi : R_1 \to R_2$ is a homomorphism of $\Ql$-groups such that
$$
\begin{CD}
\G_1 @>{f}>> \G_2 \cr
@V{\rho_1}VV @VV{\rho_2}V \cr
R_1(\Ql) @>{\phi}>> R_2(\Ql) \cr
@A{w_1}AA @AA{w_2}A \cr
\Gm(\Ql) @= \Gm(\Ql) \cr
\end{CD}
$$
commutes.

The following result is easily proved using the universal mapping property
of weighted completion and Theorem~\ref{naturality}.

\begin{proposition}
Under these hypotheses, there is a homomorphism $\Phi : \cG_1 \to \cG_2$
of proalgebraic $\Ql$-groups such that the diagram
$$
\begin{CD}
\G_1 @>f>> \G_2 \cr
@V{\rhotilde_1}VV @VV{\rhotilde_2}V \cr
\cG_1(\Ql) @>>{\Phi}> \cG_2(\Ql) \cr
\end{CD}
$$
commutes. Moreover, the induced homomorphism $\g_1 \to \g_2$ on Lie
algebras preserves the weight filtration and is strictly compatible with
it. \qed
\end{proposition}

We are now in the situation described in Paragraph~\ref{homom} and can
thus apply Theorem~\ref{naturality}.

\subsection{Right exactness}

Our goal in this paragraph is to understand the relation between
the kernel of $\cG\to R$ and the $\l$-adic unipotent completion of
the kernel of $\rho : \G \to R(\Ql)$.

We continue with the notation from the beginning of this section:
$\G$ is a profinite group, $\rho : \G \to R(\Ql)$ is a continuous, Zariski
dense homomorphism from $\G$ to an algebraic $\Ql$-group, $w$ is a
central cocharacter of $R$, $\cG$ is the corresponding weighted completion
of $\G$, and $\U$ is the kernel of $\cG \to R$.

Denote by $\cR$ the weighted completion of $\im \rho$ with respect to the
inclusion $\im\rho \hookrightarrow R(\Ql)$ and $w$. The
homomorphisms $\G \to \im \rho$ and $\ker \rho \to \U(\Ql)$ induce
homomorphisms
$$
\cG \to \cR \text{ and } (\ker\rho)^\un_{/\Ql} \to \U
$$
where $(\ker\rho)^\un_{/\Ql}$ denotes the $\l$-adic unipotent completion
(Paragraph~\ref{ladic_compln}) of $\ker \rho$.

\begin{proposition}\label{ex_seq:weighted}
The sequence
$$
(\ker\rho)^\un_{/\Ql} \to \cG \to \cR \to 1
$$
is exact.
\end{proposition}

\begin{proof}
Since $\G \rightarrow \cR(\Ql)$ is Zariski dense, $\cG \to \cR$ is surjective.
Denote the image of $(\ker\rho)^\un_{/\Ql}$ in $\U$ by $\K$. Since the
diagram
$$
\begin{CD}
\G @>>> \im \rho \cr
@V{\rho}VV @VVV \cr
\cG(\Ql) @>>> \cR(\Ql) \cr
\end{CD}
$$
commutes, the image of $\ker \rho$ in $\cR(\Ql)$ is trivial. It follows that
$\K$ is contained in the kernel of $\cG \to \cR$ and that the composite
$$
(\ker\rho)^\un_{/\Ql} \to \cG \to \cR
$$
is trivial. It remains to show that $\K$ contains the kernel of $\cG \to \cR$.

Note that $\K$ is the Zariski closure of the image of $\ker \rho$ in $\U$.
Since $\ker \rho$ is normal in $\G$, and since $\G$ is Zariski dense in $\cG$,
$\K$ is normal in $\cG$. The homomorphism $\rho$ induces a homomorphism
$\im \rho  \to (\cG/\K)(\Ql)$. Because $\cG/\K$ is a negatively weighted
extension of $R$, the universal mapping property of $\cR$ gives a splitting
$\cR \to \cG/\K$ of $\cG \to \cR$. Since the image of $\G$ is Zariski dense in
$\cG$, $\im \rho$ has Zariski dense image in $\cG/\K$, which implies that
this splitting is surjective, and therefore an isomorphism $\cR \to \cG/\K$.
\end{proof}

\begin{corollary}
\label{uni_kernel}
If $\Hcts^1(\im\rho, V_\alpha)$ vanishes for all irreducible representations
$V_\alpha$ of $R$ for which $n(\alpha)<0$, then the sequence
$$
(\ker\rho)^\un_{/\Ql} \to \cG \to R \to 1
$$
is exact. In particular, the image of $\ker \rho$ in $\U(\Ql)$ is Zariski
dense.
\end{corollary}

\begin{proof}
Corollary~\ref{ident} implies that $R=\cR$. The result follows.
\end{proof}

\section{Proof of Theorem~\ref{h1_new}}
\label{proof:h1_new}

The proof is an exercise in group cohomology. We continue with the notation
from Section~\ref{wtd_cmp}, except that we will abuse notation and not
distinguish between an algebraic $\Ql$-group and the group of its $\Ql$
rational points.

Recall that
$$
\U = \limproj \rhohat U_\rhohat
$$
where $\rhohat$ ranges over all continuous, Zariski dense representations
$\rhohat : \G \to G_\rhohat(\Ql)$ into a negatively weighted extension
$G_\rhohat$ of $R$ by a unipotent group, which we shall denote by $U_\rhohat$. 

\begin{proposition}
\label{h1:lim}
There are natural isomorphisms
$$
H_1(\u) = H_1(\U) = \limproj \rhohat H_1(U_\rhohat)
\text{ and } \Hcts^1(\u) = \liminj \rhohat H^1(\u_\rhohat).
$$
\end{proposition}

Each finite dimensional $R$-module will be considered as a continuous
$\G$-module via $\rho$. Note that there are natural one-to-one
correspondences between:
\begin{enumerate}
\item the set of continuous 1-cocycles $f : \G \to V$;
\item the set of continuous homomorphisms $\phi : \G \to R \ltimes V$ that
lift $\rho$;
\item the set of continuous splittings $s : \G \to \G\ltimes V$ of the
extension
$$
0 \to V \to \G \ltimes V \to \G \to 1.
$$
\end{enumerate}
Under this correspondence, cocycle $f$, the lift
$\gamma \mapsto (\rho(\gamma),f(\gamma))$, and the splitting
$\gamma \mapsto (\gamma,f(\gamma))$ all correspond.

The following result is somewhat standard and is easily proved. It is
a basic tool in the proof of Theorem~\ref{h1_new}.

\begin{lemma}
There are natural one-to-one correspondences between the following three sets:
$$
\left\{\parbox{1.7in}{
continuous homomorphisms $\phi : \G \to R\ltimes V$ that lift $\rho$
}\right\}\bigg/ \parbox{0.9in}{conjugation by elements of $V$}
$$
$$
\left\{\parbox{1.5in}{
continuous splittings $s : \G \to \G\ltimes V$ of $\G\ltimes V \to \G$
}\right\}\bigg/ \parbox{0.9in}{conjugation by elements of $V$}
$$
and $\Hcts^1(\G,V)$. The bijections are induced by the bijections above.
\qed
\end{lemma}

We now prove the Theorem.
First, if $V$ is a pure representation of $R$ with $n(V) \ge 0$, then
$\Homcts_S(H_1(\U),V) = 0$ as $H_1(\U)$ is negatively weighted. So suppose
that $V$ is a pure representation of negative weight. In this case, the
extension $S\ltimes V$ is negatively weighted.

Each element of $\Hcts^1(\G,V)$ gives a homomorphism
$\rhohat :\G \to R\ltimes V$ that lifts $\rho$. By the universal mapping
property of $\cG$, there is a homomorphism $\theta : \cG \to R\ltimes V$
whose composition with the canonical mapping $\G \to \cG(\Ql)$ is $\rhohat$.
Since $R\ltimes V$ is an algebraic group, this homomorphism induces a
continuous mapping $H_1(\U) \to V$, which is clearly
$R$-invariant. This homomorphism depends only on the homomorphism
$\rhohat$ up to conjugation by elements of $V$. There is thus a well
defined homomorphism
$$
\Hcts^1(\G,V) \to \Homcts_R(H_1(\U),V).
$$

Conversely, suppose that $\phi : H_1(\U) \to V$ is continuous and
$R$-invariant. Pushing out the extension
$$
0 \to H_1(\U) \to \cG/[\U,\U] \to R \to 1
$$
along $\phi$, we obtain an extension
$$
0 \to V \to G_\phi \to R \to 1
$$
which is algebraic as $\phi$ is continuous. Composing the canonical
mapping $\G \to \cG$ with $\cG \to \cG/[\U,\U] \to G_\phi$, we obtain a
homomorphism $\G \to G_\phi$ which lifts $\rho$. By Lemma~\ref{wtd_levi},
the projection $G_\phi \to R$ splits, and the splitting is
unique up to conjugation by an element of $V$. Thus $\phi$ induces a
homomorphism $\G \to R\ltimes V$ which is unique up to conjugation by an
element of $V$, and therefore an element of $\Hcts^1(\G,V)$. This
mapping is easily seen to be the inverse of the one above.
This completes the proof of Theorem~\ref{h1}.

\section{Proof of Theorem~\ref{h2}}
\label{proof:h2}

\subsection{A technical lemma}
Suppose that $R$ is a reductive group over a field $k$ of characteristic
zero. Suppose that $\n$ is a nilpotent Lie algebra over $k$ and that
$R$ acts on $\n$ as a group of automorphisms.

Suppose that $V$ is a finite dimensional $R$-module. Each element of
$$
\theta \in \Hcts^2(\n)\otimes V
$$
determines a central extension
$$
0 \to V \to \ntilde_\theta \to \n \to 0.
$$

\begin{lemma}
\label{inv_lift}
The action of $R$ on $\n$ and $V$ lifts to an action on $\ntilde_\theta$ as
a group of Lie algebra automorphisms if and only if
$$
\theta \in \big[\Hcts^2(\n)\otimes V\big]^G
$$
\end{lemma}

\begin{proof}
We will prove the sufficiency of the condition; necessity is left as an
exercise.

Since $R$ is reductive and acts on the finite dimensional complex
$$
C^\dot(\n) := (\Lambda^\dot \n^\ast, -[\blank,\blank]^\ast)
$$
of Lie algebra cochains, 
$$
\big[\Hcts^\dot(\n)\otimes V\big]^R = H^\dot(\big[C^\dot(\n)\otimes V\big]^R).
$$
In particular, each
$$
\theta \in \big[\Hcts^2(\n)\otimes V\big]^R
$$
can be represented by a continuous $R$-invariant cocycle
$\Theta : \Lambda^2 \n \to V$.
The Lie algebra $\ntilde_\theta$ can be constructed as $\n\oplus V$
with bracket
$$
[(x,v),(y,w)] = ([x,y],\Theta(x\wedge y)).
$$
The obvious action of $R$ on $\n\oplus V$ preserves the bracket as $\Theta$
is $R$ invariant.
\end{proof}

\subsection{Proof of the theorem}
We continue with the notation from Section~\ref{wtd_cmp}.
Suppose in addition that $R$ is reductive. Let $\{V_\alpha\}$ be a set of
representatives of the isomorphism classes of the finite dimensional,
irreducible $R$-modules. It suffices to construct an injective linear map
$$
\Phi_\alpha : \left[\Hcts^2(\u)\otimes V_\alpha\right]^R 
\to 
\Hcts^2(\G,V_\alpha).
$$

Suppose that $\theta \in [\Hcts^2(\u)\otimes V_\alpha]^R$. Then there exists
an algebraic quotient $G$ of $\cG$ which is an extension
$$
1 \to N \to G \to R \to 1
$$
where $N$ is unipotent, and a class
$\thetatilde \in [H^2(\n)\otimes V_\alpha]^R$, where $\n$ is the Lie algebra
of $N$, whose image under $H^2(\n) \to \Hcts^2(\u)$ is $\theta$. The class
$\thetatilde$ determines a central extension
$$
0 \to V_\alpha \to \Ntilde_\theta \to N \to 1.
$$

Choose a Levi decomposition
$$
G \cong R \ltimes N.
$$
This determines an action of $R$ on $\n$. Lemma~\ref{inv_lift} implies
that the $R$ action on $N$ lifts to $\Ntilde_\theta$. We can thus form
the extension
$$
0 \to V_\alpha \to R \ltimes \Ntilde_\theta \to R\ltimes N \to 1.
$$
Pulling this extension back along the quotient mapping
$\cG \to G = R\ltimes N$, we obtain an extension
\begin{equation}
\label{ext3}
0 \to V_\alpha \to \cGtilde_\theta \to \cG \to 1.
\end{equation}
of proalgebraic groups. Pulling this  extension back along $\G \to \cG$,
we obtain an extension
\begin{equation}
\label{ext4}
0 \to V_\alpha \to \tilde{\G}_\theta \to \G \to 1.
\end{equation}
This gives a class in $\Hcts^2(\G,V_\alpha)$.  This class is easily seen to
depend only on the class of $\theta$. This procedure therefore defines a map
$$
\Phi_\alpha :\left[\Hcts^2(\u)\otimes V_\alpha\right]^R \to 
\Hcts^2(\G,V_\alpha).
$$

If $\Phi_\alpha(\theta)$ vanishes, the extension (\ref{ext4}) splits. The
universal mapping property of $\G \to \cG(\Ql)$ then implies that (\ref{ext3})
splits. This restricts to give a splitting of the extension
$$
0 \to V_\alpha \to \tilde{\u}_\theta \to \u \to 0,
$$
where $\tilde{\u}_\theta$ denotes the Lie algebra of the unipotent radical
of $\tilde{\cG}_\theta$. It follows that $\theta$ vanishes.

\section{The Weighted Completion of $\pi_1(\Spec \OFS)$}
\label{computn}

Fix a prime number $\l$. Suppose that $F$ is a number field. Denote its
ring of integers by $\O_F$. Suppose that $S$ is a finite set of closed
points of $\Spec \O_F$ containing all primes over $\l$. Let $\OFS$ be the
set of $S$-integers of $F$, so that $\Spec \OFS = \Spec \O_F - S$. Set
$$
\GFS = \pi_1(\Spec \OFS,\Spec \overline{F}).
$$
This is the Galois group of the maximal algebraic extension of $F$
unramified outside $S$.

One has the $\l$-adic cyclotomic character $\chi_\l : \GFS \to \Zlx$. Its
image is an open subgroup of $\Zlx$ and hence of finite index.
Composing this with the inclusion $\Zlx \hookrightarrow \Gml$, we
obtain a continuous homomorphism
$$
\rho_\l : \GFS \to \Gml
$$
with Zariski dense image. As in Example~\ref{zlx}, we take
$w: \Gm \to \Gm$ be the morphism $x \mapsto x^{-2}$.

Since $\Gm$ is reductive, we can form the weighted completion
$$
\begin{CD}
\AFS(\Q_l) @>>> \Gml \cr
@A{\rhotilde_\l}AA @AA{\rho_\l}A \cr
\GFS @>>{\chi_\l}> \Zlx
\end{CD}
$$
of $\GFS$ with respect to $\rho_\l$. Here $\AFS$ is a proalgebraic
group over $\Ql$ with reductive quotient $\Gm$. Denote the prounipotent
radical of $\AFS$ by $\KFS$ and the Lie algebras of $\AFS$ and $\KFS$
by $\aFS$ and $\kFS$, respectively.

\begin{proposition}
\label{ker_cyclo_dense}
The homomorphism
$$
\ker\{ \rho_\l : \GFS \to \Gm(\Ql)\} \to  \KFS(\Ql)
$$
induced by $\rhotilde_\l : \GFS \to \AFS(\Ql)$ has Zariski dense image.
\end{proposition}

\begin{proof}
This follows directly from Example~\ref{zlx} and Corollary~\ref{uni_kernel}
as the image of $\chi_\l : \GFS \to \Zlx$ is open.
\end{proof}

\subsection{Basic structure of $\aFS$}

The Lie algebra $\aFS$, being the Lie algebra of a weighted completion,
has a natural weight filtration. Note that since $w$ is the inverse 
of the square of the standard character, all weights are even. Thus the
weight filtration of $\aFS$ satisfies
$$
\aFS = W_0 \aFS, \quad \kFS= W_{-2}\aFS
\text{ and } \Gr^W_{2n+1} \aFS = 0 \text{ for all }n.
$$

\begin{theorem}
\label{a_l}
The Lie algebra $\Gr^W_\dot\kFS$ is a free Lie algebra and there is a natural
$\Gm$-equivariant isomorphism
$$
\Hcts^1(\kFS) \cong \bigoplus_{n=1}^\infty\Hcts^1(\GFS,\Ql(n))\otimes\Ql(-n)
\cong \bigoplus_{n=1}^\infty \Ql(-n)^{d_n},
$$
where
$$
d_n =
\begin{cases}
r_1 + r_2 + \#S - 1 & \text{ when } n = 1; \cr
r_1 + r_2 & \text{ when $n$ is odd and $> 1$;} \cr
r_2 & \text{ when $n$ is even.}
\end{cases}
$$
Here $r_1$ and $r_2$ are the number of real and complex places of $F$,
respectively. Moreover, each weight graded quotient of $\aFS$ is a finite
dimensional $\Ql$-vector space.
Any lift of a graded basis of $H_1(\Gr^W_\dot \kFS)$ to a graded
set of elements of $\Gr^W_\dot \kFS$ freely generates $\Gr^W_\dot \kFS$.
\end{theorem}

\begin{proof}
This follows directly from Theorem~\ref{h1}, Corollary~\ref{freeness},
and Theorem~\ref{Th:Galcoh}.
\end{proof}

\subsection{The case $F=\Q$}
Set
$$
G_\l = G_{\Q,\{\l\}},\, \al = \a_{\Q,\{\l\}},\, \kl = \k_{\Q,\{\l\}},\,
\Al = \A_{\Q,\{\l\}},\, \Kl = \K_{\Q,\{\l\}}.
$$

Combining the results of the previous paragraph, we have:

\begin{theorem}
\label{q_gens}
The Lie algebra $\kl$ is free and
$$
\Hcts^1(\kl) = \bigoplus_{n \ge 0} \Ql(-2n-1).
$$
Moreover, there are elements $\sigma_m \in W_{-2m}\al$ for each odd positive
number $m$ whose images in $\Gr^W_\dot \kl$ generate it freely. \qed
\end{theorem}

\subsection{The weight filtration of $\GFS$}

The weight filtration of $\aFS$ can be exponentiated to give a weight
filtration
$$
\AFS = W_0 \AFS \supseteq W_{-1} \AFS \supseteq W_{-2} \AFS
\supseteq W_{-3}\AFS \supseteq \cdots
$$
of $\AFS$ whose odd weight graded quotients vanish. This induces a filtration
on $\GFS$ by pulling back the weight filtration
$$
\GFS = W_0 \GFS \supseteq W_{-1} \GFS = W_{-2} \GFS
\supseteq W_{-3}\GFS = W_{-4} \GFS \supseteq \cdots
$$
of $\AFS$ along the natural map $\rhotilde_\l :\GFS \to \AFS(\Ql)$. Note that
$W_{-1}\GFS$ is the kernel of the cyclotomic character.

For each $m>0$, there is a natural continuous inclusion
\begin{equation}
\label{inclusion}
\Gr^W_{-m} \GFS \hookrightarrow \Gr^W_{-m} \AFS(\Ql)
\stackrel{\sim}{\longrightarrow}\Gr^W_{-m} \aFS.
\end{equation}

\begin{theorem}
\label{Th:GrGl-Al}
The linear inclusions (\ref{inclusion}) combine to give an isomorphism
$$
\Gr^W_{<0}\GFS \otimes_{\Zl}\Ql \hookrightarrow \Gr^W_\dot \kFS
$$
of $\Ql$-Lie algebras.
\end{theorem}

Since, by Proposition~\ref{ker_cyclo_dense}, the image of $\ker \chi_\l$ in
$\KFS(\Ql)$ is Zariski dense, this follows from the general lemma proved in
the next paragraph.

\subsection{Another technical lemma}

Suppose that $\G$ is a profinite group, that $U$ is a prounipotent
$\Ql$-group and that
$$
U = W_{-1}U \supseteq W_{-2}U \supseteq W_{-3}U \supseteq \cdots
$$
is a central filtration of $U$ such that $U/W_{-n}U$ is unipotent for each
$n$. Suppose that $\rho : \G \to U(\Ql)$ is a continuous homomorphism. We
can induce a filtration
$$
\G = W_{-1}\G \supseteq W_{-2}\G \supseteq W_{-3}\G \supseteq \cdots
$$
on $\G$ by defining $W_{-n}\G = \rho^{-1}(W_{-n}U(\Ql))$. Note that the
mapping
\begin{equation}
\label{incln}
\Gr^W_m \G \to \Gr^W_m U(\Ql)
\end{equation}
is injective.

\begin{lemma}
\label{graded_iso}
If the image of $\rho:\G \rightarrow U(\Ql)$ is Zariski dense, then the
inclusion (\ref{incln}) induces an isomorphism
$$
\Gr_{m}^W\G\otimes_{\Zl}\Ql 
\rightarrow \Gr^W_{m} U(\Ql)
$$
for each $m <  0$.
\end{lemma}

\begin{proof}
Fix $m < 0$.
By replacing $U$ by $U/W_M U$ for some $M < m$, we may assume that $U$ is
unipotent. By replacing $\G$ by $\G/W_M\G$, we may assume that $\rho$
is an inclusion.

We will prove the result by induction on the dimension of $U$. First recall
that if $U$ is a unipotent group over a field $k$ of characteristic zero,
and $V$ is a closed normal subgroup of $U$, also defined over $k$, then the 
group of $k$-rational points of $U/V$ is isomorphic to $U(k)/V(k)$.

We first consider the case where $U$ is abelian. In this case, $\G$ is a
compact subgroup of $U(\Ql) \cong \Ql^N$. Since $\Zl$ is a PID, and since
$\G$ is compact and torsion free, it is freely generated by $r$ linearly
independent elements of $\Ql^N$. It follows that $\G$ is Zariski dense
in $\Ql^N$ if and only if $r=N$. This proves the result when $U$ is abelian.

Now suppose that $U$ is not abelian. Note that the commutator subgroup
$[\G,\G]$ is Zariski dense in $[U,U]$; it is clear that the Zariski closure
of $[\G,\G]$ is contained in $[U,U]$. The reverse conclusion follows as
the image of $\G$ in the $\Ql$-rational points of the quotient of $U$
by the closure of $[\G,\G]$ is dense and abelian, which implies that
$$
[U,U] \subseteq \text{ the Zariski closure of } [\G,\G] \text{ in } U.
$$

If $U$ is not abelian, then there is a least $n$ for which $U \neq W_{-n} U$.
Since the filtration $W_\dot$ is central, $W_{-n} U \supseteq [U,U]$. Using
the fact that $[\G,\G]$ is Zariski dense in $[U,U]$ and also the fact that
the result holds for $U/[U,U]$, it is not hard to see that $W_{-n}\G$ is
Zariski dense in $W_{-n}U$. By induction, the result holds for
$$
W_{-n}\G \hookrightarrow W_{-n}U(\Ql)
$$
from which the result for $\G \to U(\Ql)$ follows as
$$
\Gr^W_{-m} U =
\begin{cases}
\Gr^W_{-m} U/W_{-n} U & m < n, \cr
\Gr^W_{-m} W_{-n} U  & m \ge n.
\end{cases}
$$
\end{proof}

\section{Galois Actions on Geometric Fundamental Groups}
\label{actions}

Fix a prime number $\l$. Suppose that $F$ is a number field and
that $X$ is a geometrically connected variety defined over $F$. Set
$\Xbar = X\otimes_F \Fbar$. Denote the
absolute Galois group of $F$ by $G_F$. We have the outer Galois action
$\phi_\l : G_F \to \Out \pi_1(\Xbar,x)^\prol$
of $G_F$ on the pro-$\l$ completion of the geometric fundamental group
of $(X,x)$.

For a finite set of primes $S$ of the ring of integers $\O_F$,
we denote the ring of $S$-integers by $\OFS$,
and $\pi_1(\Spec \OFS, \Spec\overline{F})$ by $\GFS$.
It is known that there is a finite set of primes $S$
such that $\phi_\l$ factors through $G_F \to \GFS$, so
we have
$$
\phi_S : \GFS \to \Out \pi_1(\Xbar,x)^\prol.
$$
This $S$ can be taken as the union of  the bad-reduction primes of $X$
and the primes above $\l$, (cf.\ \cite[Chapter XIII]{SGA1}). If $X=\Pminus$
over $F=\Q$, then we can take $S=\{\l\}$.

\subsection{The $I$-filtration}
Define a filtration 
$$
\Out \pi_1(\Xbar,x)^\prol =
L^0 \Out\pi_1(\Xbar,x)^\prol \supseteq L^1 \Out \pi_1(\Xbar,x)^\prol
\supseteq \cdots
$$
of the outer automorphism group of $\pi_1(\Xbar,x)^\prol$ by
\begin{multline*}
L^m \Out \pi_1(\Xbar,x)^\prol \\
= \ker\big\{\Out \pi_1(\Xbar,x)^\prol \to
\Out \big(\pi_1(\Xbar,x)^\prol/L^{m+1}\pi_1(\Xbar,x)^\prol\big)\big\}
\end{multline*}
where
$$
\pi_1(\Xbar,x)^\prol = L^1\pi_1(\Xbar,x)^\prol \supseteq
L^2\pi_1(\Xbar,x)^\prol \supseteq \cdots
$$
denotes the lower central series of $\pi_1(\Xbar,x)^\prol$.
The $I$-filtration
$$
G_F = \I^0 G_F \supseteq \I^1 G_F \supseteq \I^2 G_F \supseteq
\cdots
$$
of $G_F$ is the pull back of the filtration
$L^\dot$ of $\Out \pi_1(\Xbar,x)^\prol$ along the outer action:
$$
\I^m G_F = \phi_\l^{-1} L^m \Out \pi_1(\Xbar,x)^\prol
= \ker \big\{G_F \to \Out\big(\pi_1(\Xbar,x)^\prol/L^{m+1}\big)\big\}.
$$
We can also pullback the filtration $L^\dot$ of $\Out \pi_1(\Xbar,x)^\prol$
to define the $I$-filtration of $\GFS$:
$$
\GFS = \I^0 \GFS \supseteq \I^1 \GFS \supseteq \I^2 \GFS 
\supseteq \cdots
$$
By an elementary argument, the quotient mapping
$G_F \twoheadrightarrow \GFS$ induces isomorphisms
\begin{equation}
\label{gr_isom}
\Gr_\I^n \GFS \cong \Gr_\I^n G_F.
\end{equation}
Since $[\I^m G_F,\I^n G_F]\subseteq \I^{m+n}G_F$,
the associated graded groups
$$
\Gr_\I^{>0} G_F \cong \Gr_\I^{>0} \GFS
$$
are Lie algebras over $\Zl$; the bracket is induced by the group
commutator. (See \cite{serre:LALG}, for example.)

\subsection{Galois actions on unipotent completions}

Set $\cP = \pi_1(\Xbar,x)^\un_{/\Ql}$, the $\l$-adic unipotent completion
(see Paragraph~\ref{ladic_compln}) of $\pi_1(\Xbar,x)$.

\begin{remark}
It follows from Theorem~\ref{comp_thm} that for each imbedding
$\sigma : F \hookrightarrow \C$, there is an isomorphism
$$
\cP = \pi_1(X_\sigma,x)^\un_{/\Q}\otimes_\Q \Ql,
$$
where $X_\sigma$ denotes the complex variety obtained from $X$ via
$\sigma$.
\end{remark}

There is a proalgebraic group $\cOut \cP$ whose $K$-rational points form
the group of continuous outer automorphisms of $\cP(K)$ for each field
extension $K$ of $\Ql$. (See Paragraph~\ref{autos}.) By the functoriality of
unipotent completion, the outer action of $\GFS$ on $\pi_1(\Xbar,x)$
induces a homomorphism $\Phi_S : \GFS \to \cOut \cP$.
As above, the lower central series
$$
\cP = L^1 \cP \supseteq L^2 \cP \supseteq L^3 \cP \supseteq \cdots
$$
of $\cP$ induces a filtration
$$
L^0 \cOut \cP \supseteq L^1\cOut \cP \supseteq L^2\cOut \cP \supseteq \cdots
$$
of $\cOut \cP$. Pulling this back along $\Phi_S$, we obtain a filtration
$$
\GFS = \cI^0 \GFS \supseteq \cI^1 \GFS \supseteq
\cI^2 \GFS \supseteq \cdots
$$
of $\GFS$.

\begin{proposition}\label{prop:comp-gr}
For all $m \ge 0$, $\I^m \GFS \subseteq \cI^m \GFS$ and the natural map
$$
\Gr_\I^m \GFS \to \Gr_\cI^m \GFS
$$
has finite kernel and cokernel. In particular,
if $m>0$, then
$$
(\Gr_\I^m \GFS) \otimes_\Zl \Ql \cong (\Gr_\cI^m \GFS) \otimes_\Zl \Ql.
$$
\end{proposition}

\begin{proof}
We will use Theorem~\ref{ladic_lcs_gradeds} and its notation. Set
$\G = \pi_1(\Xbar,x)^\prol$. Since $\G$ is a finitely generated pro-$\l$
group, $L^m \G \subseteq D^m\G$ and $D^m\G/L^m\G$ is a finite $\l$-group.
Also, $L^m\G = D^m\G$ for all $m$ if and only if $\G/L^m\G$ is torsion
free for all $m$.

Using universal mapping properties, one can show that the inclusion
$\G/D^m\G \hookrightarrow (\cP/L^m\cP)(\Ql)$ is the $\l$-adic
unipotent completion of $\G/D^m\G$. It follows that the
natural homomorphism $\Aut \G \to \Aut \cP$ descends to a natural
homomorphism
$$
\Aut(\G/D^m\G) \hookrightarrow \Aut(\cP/L^m\cP)
$$
which is injective. This induces a homomorphism
$$
\Out(\G/D^m\G) \to \Out(\cP/L^m\cP).
$$
The kernel of this morphism is finite. This is equivalent to the assertion
that $\Inn(\G/D^m\G)$ has finite index in
$\Aut(\G/D^m\G) \cap \Inn(\cP/L^m\cP)$, which follows from the compactness
of $\Aut(\G/D^m\G)$ ($\G$ is finitely generated as a pro-$\l$ group) and the
openness of $\G/D^m\G$ in $\cP/L^m\cP(\Ql)$, a consequence of
Lemma~\ref{graded_iso}. Since $D^m\G/L^m\G$ is a finite $\l$-group, it
follows that the kernel of both
$$
\Out(\G/L^m\G) \to \Out(\G/D^m\G) \text{ and }
\Out(\G/L^m\G) \to \Out(\cP/L^m\cP)
$$
are finite from which it follows that
$\I^m \GFS$ is a finite index subgroup of $\cI^m \GFS$.

Because the sequence
$$
1 \to (\I^m \cap \cI^{m+1})/\I^{m+1} \to \I^m/\I^{m+1} 
\stackrel{\phi}{\to} \cI^m/\cI^{m+1} \to \cI^m/(\I^m + \cI^{m+1}) \to 1
$$
is exact, the kernel and cokernel of $\phi$ are finite.

Since $[\I^m,\I^n] \subseteq \I^{m+n}$ and $[\cI^m,\cI^n] \subseteq \cI^{m+n}$,
$\Gr_\I^m \GFS$ and $\Gr_\cI^m \GFS$ are abelian pro-$\l$ groups whenever
$m>0$. Thus it follows that the kernel and cokernel of the mapping on each
associated graded is a finite abelian $\l$-group when $m>0$.
\end{proof}

\subsection{The Galois image}
Choose an $F$-rational base point, --- either a geometric point
$x$ or a tangential base point anchored at an $F$-rational
point of a smooth completion of $X$. Denote the Lie algebra of $\cP$ by $\p$. 
It is well known that $\Het^1(\Xbar,\Ql)$ is isomorphic to 
$\Homcts(\pi_1(\Xbar,x),\Ql)$ as a continuous $\GFS$-module.

Denote the Zariski closure of the image of $\GFS$ in $\cAut\cP$ by
$\cGhat_{F,S}$ and its prounipotent radical by $\Uhat_{F,S}$. Denote the
Zariski closure of the image of $\GFS$ in $\cOut\cP$ by $\cG_{F,S}$ and
its prounipotent radical by $\U_{F,S}$. There is a natural homomorphism
$\cGhat_{F,S} \to \cG_{F,S}$.

\begin{proposition}
\label{action}
If $\Het^1(\Xbar,\Ql(1))$ is a trivial $\GFS$-module, then
\begin{enumerate}
\item the natural homomorphism $\GFS \to \cAut \cP \to \Aut H_1(\p)$
is the composition of the $\l$-adic cyclotomic character with the inclusion
of the scalar matrices;
\item the subgroup $\cGhat_{F,S}$ of $\cAut\cP$ is a negatively weighted
extension of the scalar matrices $\Gm$ by $\Uhat_{F,S}$ with respect to 
the central cocharacter $w : x \mapsto x^{-2}$;
\item the subgroup $\cG_{F,S}$ of $\cOut\cP$ is a negatively weighted
extension of the scalar matrices $\Gm$ by $\U_{F,S}$ with respect to 
the central cocharacter $w : x \mapsto x^{-2}$;
\item The homomorphism $\GFS \to \Aut\U(\Ql)$ induces surjective homomorphisms
$$
\AFS \to \cGhat_{F,S} \to \cG_{F,S};
$$
\item $\p$ has a natural weight filtration which is essentially its lower
central series:
$$
W_{-2m} \p = L^m \p \text{ and } Gr^W_{2m+1}\p = 0 \text{ for all } m;
$$
\item\label{der} the weight filtrations on $\p$, $\Der\p$ and $\OutDer \p$
are related by
$$
\Gr^W_n \OutDer \p = \left(\Gr^W_n \Der \p\right)/\left(\Gr^W_n \p\right)
\text{ and }
\Gr^W_\dot \Der\p \cong \Der \Gr^W_\dot \p;
$$
\item the weight filtration of $\Uhat_{K,T}$ is characterized by
$$
W_{-2m}\, \Uhat_{F,S} =
\ker\big\{ \Uhat_{F,S} \to \cAut (\cP/L^{m+1}\cP) \big\}
\text{ and } \Gr^W_{2m+1}\Uhat_{F,S} = 1
$$
for all $m$;
\item \label{defn}
the weight filtration of $\U_{F,S}$ is characterized by
$$
W_{-2m}\, \U_{F,S} = \ker\big\{ \U_{F,S} \to \cOut (\cP/L^{m+1}\cP) \big\}
\text{ and }
\Gr^W_{2m+1}\U_{F,S} = 1
$$
for all $m$.
\end{enumerate}
\end{proposition}

\begin{proof}
The first assertion follows from the assumption using the natural
isomorphism $H_1(\Xbar) \cong H_1(\p)$. Since the bracket
$$
\Gr_L^1 \p \otimes \Gr_L^n \p \to \Gr_L^{n+1} \p
$$
is surjective and $\GFS$ equivariant, it follows that $\GFS$ acts on
$\Gr_L^m \p$ via the $m$th power of the cyclotomic character. This implies
the second and third assertions as it implies that every derivation of
$\p$ that acts trivially on $\Gr^\dot_L \p$ has negative weight. The fourth
follows from the universality of $\AFS$.
Because $\p$ is a module over $\AFS$, it has a natural weight filtration.
Since $\GFS$ acts on $\Gr_L^n \p$ via the $n$th power of the cyclotomic
character, next assertion follows. The remaining assertions follow as
$\Der \p$ and $\OutDer \p$ are negatively weighted $\GFS$-modules.
\end{proof}

\begin{theorem}
\label{generation}
If $\Het^1(\Xbar,\Ql(1))$ is a trivial $\GFS$-module, then there
is a graded Lie algebra surjection
$\Gr^W_\dot \k_{F,S} \to \Gr^W_\dot \u_{F,S}$ and an isomorphism
$$
\Gr^W_\dot \u_{F,S} \cong (\Gr_\I^{>0} \GFS)\otimes_\Zl \Ql.
$$
\end{theorem}

\begin{proof}
The surjection $\A_{F,S} \to \cG_{F,S}$ restricts to a surjection
$\KFS \to \U_{F,S}$, which implies that $\I^1\GFS \to \U_{F,S}(\Ql)$ has
Zariski dense image. By strictness, it
induces a surjection $\Gr^W_\dot \k_{F,S} \to \Gr^W_\dot \u_{F,S}$.
{}From Proposition~\ref{action}, it follows that $\cI^m \GFS$
is the inverse image of $W_{-2m} \cG_{F,S}$ under the natural mapping
$\GFS \to \cG_{F,S}$. There is thus an injection
$$
\Gr^m_\cI \GFS \hookrightarrow \Gr^W_{-2m} \cG_{F,S}.
$$
Proposition~\ref{prop:comp-gr}, Proposition~\ref{ker_cyclo_dense}
and Lemma~\ref{graded_iso} imply
that this induces isomorphisms
$$
(\Gr^m_\I \GFS)\otimes_\Zl \Ql \cong
(\Gr^m_\cI \GFS)\otimes_\Zl \Ql \cong
\Gr^W_{-2m} \U_{F,S} \cong \Gr^W_{-2m} \u_{F,S}.
$$
whenever $m>0$. 
\end{proof}

\subsection{Proof of Conjecture~\ref{main}}
Since $\Pminus$ has everywhere good reduction, the Galois representation
$$
\phi_\l : G_\Q \to \Out \pi_1(\P^1_{\Qbar} - \{0,1,\infty\})^\prol
$$ 
factors through the projection 
$G_\Q \to G_\l=G_{\Q, \{\l\}}$ \cite[Thm.~1]{Ihara1}.
It is standard that $\Het^1(\P^1_\Qbar - \{0,1,\infty\},\Ql(1))$
is a trivial $G_\l$-module. Now apply Theorem~\ref{generation} with
$X=\Pminus$, $F=\Q$, $S=\{\l\}$, $\kFS=\kl$. Theorem~\ref{q_gens} says that
$\Gr^W_\dot \kl$ is generated by elements
$$
\sigma_{2n+1} \in \Gr^W_{-2(2n+1)} \kl \quad n \ge 0.
$$
Conjecture~\ref{main} now follows from Theorem~\ref{generation} and the 
following lemma, which implies that the image of $\sigma_1$ in
$$
(\Gr_I^{>0} \Gl)\otimes\Ql \cong \Gr^W_\dot \u_{\Q,\{\l\}}
$$
is trivial.

\begin{lemma}
\label{lem:sigma1}
We have $\Gr^W_{-2} \OutDer \p = 0$, so that $\Gr^W_{-2}\u_{\Q,\{\l\}} = 0$.
\end{lemma}

\begin{proof}
Since $\Gr^W_\dot \p$ is isomorphic to the free Lie algebra generated
by $H_1(\p)$, it follows that $\Gr^W_{-4} \p \cong \Ql(2)$. By
Proposition~\ref{action} we have
$$
\Gr^W_{-2}\Der \p \cong \Gr^W_{-2}\Der\Gr^W_\dot\p \cong
\Hom(H_1(\p),\Gr^W_{-4}\p) \cong H_1(\p)
$$
The non-degeneracy of the Lie bracket 
$$
H_1(\p) \otimes H_1(\p) \to \Gr^W_{-4} \p \cong \Ql(2)
$$
implies that the mapping $H_1(\p) \to \Hom(H_1(\p), \Gr^W_{-4}\p)$, given
by taking inner derivations, is an isomorphism. Thus
$$
\Gr^W_{-2}\OutDer \p \cong \Gr^W_{-2}\Der \p/\Gr^W_{-2}\p 
\cong \Hom(H_1(\p), \Gr^W_{-4}\p)/H_1(\p)=0.
$$
\end{proof}

\begin{remark}
\label{h1_compn}
Our proof says nothing about the non-triviality of the
$$
\sigma_{2m+1} \in (\Gr_\I^\dot \Gl)\otimes\Ql.
$$
Ihara's work \cite{Ihara3} implies that when $m \ge 1$, all are non-trivial
in $H_1((\Gr_\I^\dot \Gl)\otimes\Ql)$. Combined with Theorem~\ref{q_gens},
this implies that
$$
H_1((\Gr_\I^\dot \Gl)\otimes\Ql) = \bigoplus_{m\ge 1} \Ql(2m+1)
$$
where the copy of $\Ql(2m+1)$ is spanned by $\sigma_{2m+1}$.
Ihara \cite{Ihara3} also uses power series methods to
establish nonvanishing results for some brackets of the $\sigma_j$.
Improvements can be found in \cite{matsumoto} and \cite{tsunogai}.

Ihara defined a certain $\Z$-lie algebra $\cD$ 
called the {\it stable derivation algebra},
which is a subalgebra of outer derivations of 
free $\Z$-Lie algebra with two generators.
It is proved that $\Gr_I^{>0} \Gl$ is contained in
$\cD\otimes \Zl$ for every prime $\l$.
He considered the problem of whether $\cD$
is generated by certain derivations $D_{2m+1}$,
which are analogues of the $\sigma_{2m+1}$.
He found a mysterious congruence \cite{Ihara9}
$$
2[D_3, D_9]-27[D_5, D_7] \equiv 0 \mod 691
$$
when $\l=691$. This suggests that one might have
$$
2[\sigma_3, \sigma_9]-27[\sigma_5, \sigma_7] \equiv 0 \mod 691
$$
in $\Gr_\I^{>0} \Gl$ for $\l=691$, hence
this Lie algebra might {\it not} be generated 
by the elements $\sigma_i$ over $\Zl$ for $\l=691$. 
This is not currently known, but
Sharifi announced an interesting progress:
for $\l$ regular, the freeness of $\Gr_\I^{>0} \Gl\otimes \Ql$
implies that $\Gr_\I^{>0}$ itself is generated by 
the elements $\sigma_i$, and for $\l$ irregular, a conjecture
of Greenberg's implies that $\Gr_\I^{>0}\Gl$ is not free on $\sigma_i$.
\end{remark}

\subsection{Goncharov's Conjecture}
\label{goncharov}
Goncharov \cite{goncharov} considers the varieties
$$
X_N := \P^1 - \{0,\mu_N,\infty\}
$$
where $\mu_N$ denotes the group of $N$th roots of unity. Take $F$ to
be $\Q(\mu_N)$ and $S$ to be the set of primes in $\O_F$ that lie over
$N\l$. Let $v = \overrightarrow{01}$, the tangent vector at 0
that points towards 1. The Galois representation
$$
\phihat_\l : G_F \to \Aut \pi_1(\Xbar,v)^\prol
$$
is unramified outside $S$ by the Smooth Base Change Theorem
\cite[Chapt.~XIII]{SGA1}, so it factors through a representation
$$
\phihat_S : \GFS \to \Aut \pi_1(\Xbar,v)^\prol.
$$

Denote the $\l$-adic unipotent completion of $\pi_1(\Xbar,v)^\prol$ by
$\cP$. By functoriality, $\phihat_\l$ and $\phihat_S$ induce homomorphisms
$G_F \to \GFS \to \Aut \cP(\Ql)$.
By Proposition~\ref{action}, this induces a homomorphism
$$
\Phihat_S : \AFS \to \cAut\cP.
$$

Goncharov considers the Zariski closure of the image of
$$
G_{\Q(\mu_{\l^\infty})} := \Gal(\Qbar/\Q(\mu_{\l^\infty}))
$$
in $\cAut\cP$. Since
$$
G_{\Q(\mu_{\l^\infty})} = \ker\{\chi_\l : G_F \to \Zlx\},
$$
Theorem~\ref{ker_cyclo_dense} implies that this is the image $\Uhat_{F,S}$
of $\KFS$ in $\cAut\cP$.

\begin{theorem}
For all $N \ge 1$, the Lie algebra of the Zariski closure of the image
of $G_{\Q(\mu_{\l^\infty})} \to \cAut\cP$ is a quotient of $\kFS$, and
is therefore generated topologically by a lift of
$$
\bigoplus_{n\ge 0} \Het^1(\Spec \O_{{\Q(\mu_N)},S},\Ql(n))^\ast \otimes \Ql(n).
$$
\end{theorem}

When $N=1$, this proves the generation portion of Conjecture~2.1 in
\cite{goncharov}.

\section{$\l$-adic Mixed Tate Modules over $\Spec \OFS$}
\label{deligne}

In this section, we show how weighted completion of Galois groups can
be used to prove $\l$-adic versions of Deligne's Conjectures
\cite[Sect.~8]{deligne} on mixed Tate motives over the spectrum of the ring
of $S$-integers in a number field.
Such $\l$-adic versions of these conjectures have previously been proved by
Beilinson and Deligne by different, but equivalent, methods in their
unpublished manuscript \cite{beilinson-deligne}. 

\subsection{Deligne's conjectures}
Suppose that $F$ is a number field and that $S$ is a finite subset of
$\Spec \O_F$. Deligne \cite[8.2]{deligne} conjectures that there
is a category of mixed Tate motives over $\Spec \OFS$ and that in this
category, $\Ext^1(\Q(0),\Q(n)) = K_{2n-1}(\OFS)\otimes \Q$.
The tannakian fundamental group of this category is of the form
$\Gm \ltimes U$ where $U$ is prounipotent. He further conjectures
\cite[8.9.5]{deligne} that the weight graded Lie algebra of $U$ is
freely generated by
$$
\bigoplus_{n\ge 1} K_{2n-1}(\OFS)\otimes\Q.
$$
(The generation statement is equivalent to the first conjecture above.)

The analogue of this conjecture for mixed Tate motives over a number field
has been solved by Goncarov \cite{goncharov:mtm}, who used results of Voevodsky
\cite{voevodsky} and Levine \cite{levine} to construct the category
of mixed Tate motives over a number field. He has also made a working
definition of mixed Tate motives over the ring of $S$-integers in a number
field. 

Here we formulate and prove an $\l$-adic version of this. Fix a rational
prime $\l$. Beilinson and Deligne had given a similar construction in
the unpublished manuscript \cite{beilinson-deligne}.

\subsection{The case $\{\l\}$ contained in $S$}

First, we shall assume that $S$ contains all primes that lie over $\l$.
Later we shall remove this hypothesis. In this case, it is natural to
define the category of $\l$-adic mixed Tate modules over $\Spec\OFS$ to
be the category of weighted $\GFS$-modules 
(see Definition~\ref{def:wtG-mod}). Denote this category by $\MTMl(\OFS)$.

\begin{lemma}
The category $\MTMl(\OFS)$ is equivalent to the category of finite
dimensional representations of $\AFS$. For each $n\in \Z$ and $i \ge 0$,
there is a natural isomorphism
$$
\Ext^i_{\MTMl(\OFS)}(\Ql,\Ql(n)) \cong H^i(\AFS,\Ql(n)).
$$
\end{lemma}

\begin{proof}
This follows directly from \cite{hain-matsumoto:exp} --- see
Remark~\ref{tannaka}.
\end{proof}

If $\U$ is prounipotent with Lie algebra $\u$, then $H^\dot(\U) \cong
\Hcts^\dot(\u)$, since the category of $\U$-modules
is equivalent to the category of continuous $\u$-modules.

\begin{theorem}
\label{Th:MTMl}
There are natural isomorphisms
$$
\Ext^i_{\MTMl(\OFS)}(\Ql,\Ql(n)) \cong
\begin{cases}
\Ql & \text{ when $i = n = 0$,} \cr
\Hcts^1(\GFS,\Ql(n)) & \text{ when $i=1$ and $n>0$,} \cr
0  & \text{ otherwise.}
\end{cases}
$$
Consequently, for all $n\in \Z_{>0}$, there is a natural isomorphism
$$
\Ext^1_{\MTMl(\OFS)}(\Ql,\Ql(n)) \cong K_{2n-1}(\Spec\OFS)\otimes\Ql.
$$
\end{theorem}

\begin{proof}
Because of the previous lemma, we only need compute $H^i(\AFS,\Ql(n))$.
It is proved in \cite{jantzen} that there is a Leray-Serre type spectral
sequence
$$
H^s(\Gm,H^t(\KFS,\Ql(n))) \implies H^{s+t}(\AFS,\Ql(n)).
$$
Since $\Gm$ is reductive and $H^\dot(\KFS)\cong \Hcts^\dot(\kFS)$, it follows
that
$$
H^i(\AFS,\Ql(n))=[\Hcts^i(\kFS)\otimes \Ql(n)]^\Gm.
$$
Since $\kFS$ is negatively weighted, the right hand side is zero when
$n\leq 0$. Since $\Gr^W_\dot\kFS$ is free and $\Gr^W_\dot$ is exact,
$$
\Gr^W_\dot \Hcts^i(\kFS) = H^i(\Gr^W_\dot\kFS) = 0 
$$
when $i>1$ and $\Hcts^i(\kFS)$ vanishes if $i\geq 2$. When $i=1$,
Theorem~\ref{a_l} gives the result. For $i=0$, the assertion is obvious.

The last assertion follows from Proposition~\ref{prop:comparison} and the
fact that the regulator mappings
$$
c_1 : K_{2n-1}(\OFS)\otimes\Ql \to \Het^1(\Spec \OFS,\Ql(n))
$$
are isomorphisms for all $n$ and $\l$. This is due to Soul\'e \cite{soule}
when $\l\neq 2$ and Rognes and Weibel \cite{Rognes} when $\l = 2$.
\end{proof}

\subsection{Constrained weighted completion}
\label{sec:withproperty}

In order to handle the case where $\l$ is not contained in $S$, we
need to consider a variant of weighted completion. We use the notation
of Section~\ref{wtd_cmp}. Suppose that $P$ is a property of weighted
$\G$-modules that is closed under direct sums, tensor products, taking
duals, and taking subquotients. Suppose also that the trivial representation
$k$ has property $P$.

The category of finite dimensional weighted $\G$-modules with property $P$
is a full subcategory of the category of weighted $\G$-modules. In fact, it
is a neutral full subtannakian category of the category of weighted
$\G$-modules. The {\it weighted completion of $\G$ constrained by $P$} is
defined to be the tannakian fundamental group $\cG_P $ of this subcategory.
Denote the kernel of $\cG_P \to R$ by $\U_P$ and its Lie algebra by $\u_P$.

By construction, the categories of weighted $\G$-modules with property $P$ 
and the category of finite dimensional $\cG_P$-modules are equivalent.
Suppose that $V$ is a weighted $\G$-module with property $P$. Denote by
$\Hcts^1(\G,V)_P$ the subgroup of $\Hcts^1(\G,V)$ generated by the classes
of extensions
$$
0 \to V \to E \to k \to 0
$$
of $\G$-modules in which $E$ is also a weighted $\G$-module with property $P$.

\begin{proposition}
\label{prop:h1-P}
For all negatively weighted $\G$-modules $V$ with property $P$, there is a
natural isomorphism
$$
H^1(\cG_P,V) \stackrel{\simeq}{\longrightarrow} \Hcts^1(\G,V)_P.
$$
When $R$ is reductive and each $\Hcts^1(\G,V_\alpha)$ is finite dimensional,
$$
H^1(\U_P)\cong 
\Hcts^1(\u_P)\cong \bigoplus_{\{\alpha : n(\alpha) < 0 \}}
\Hcts^1(\G,V_\alpha)_P\otimes V_\alpha^\ast.
$$
Here we are using the notation of Paragraph~\ref{reductive}.
\end{proposition}

\begin{proof}
The first statement follows because $\Hcts^1(\G,V)_P$ classifies extensions
of $k$ by $V$ with property $P$. Since $V$ is negative, each extension of
$k$ by $V$ is a weighted $\G$-module and hence a $\cG_P$-module. It therefore
determines an element of $H^1(\cG_P,V)$.

The second assertion follows using the Leray-Serre Spectral Sequence
\cite{jantzen} for the extension
$$
1 \to \U_P \to \cG_P \to R \to 1
$$
and the fact that $R$ is reductive, which implies that
$H^s(R,H^t(\U_P,V_\alpha))$ vanishes whenever $s>0$. We
also use the fact that $\Hcts^\dot(\u_P) \cong H^\dot(\U_P)$.
\end{proof}

To obtain the analogue of Theorem~\ref{h2} for $\Hcts^2(\u_P)$, we
need to impose an extra condition on $P$.

\begin{definition}
\label{def:adjacent}
Assume that there is a natural number $b$ such that for all weighted
$\G$-modules $V$ $\Gr^W_m V=0$ when $m$ is not a multiple of $b$.
We say that property $P$ of weighted $\G$-modules has the {\it bootstrap
property with index $b$} if $V$ has property $P$ whenever $W_{m}V/W_{m-2b}V$
has property $P$ for all $m$.
\end{definition}

Note that the assumption is always satisfied when $b=1$ and that it is
satisfied when $b = 2$ in the mixed Tate case.

\begin{theorem}\label{th:h2-P}
Assume that $R$ is reductive and that $\Hcts^2(\G,V_\alpha)$ is finite
dimensional for each irreducible representation $V_\alpha$ of $R$.
Let $P$ be a property of weighted $\G$-modules that it is closed under
tensor products, direct sums, and taking duals and subquotients
as $\G$-modules. Suppose also that the trivial module $k$ has property
$P$. If there is a $b\ge 1$ such that $P$ has the bootstrap property with
index $b$, then there is a natural injection
$$
\Hcts^2(\u_P) \hookrightarrow
\bigoplus_{\{\alpha:n(\alpha)\leq -2n\}}
\Hcts^2(\G,V_\alpha)\otimes V_\alpha^\ast,
$$
where $n$ is any positive integer such that $W_{n-1}\Hcts^1(\u_P)=0$
or, equivalently, $W_{-n}\u_P = \u_P$.
\end{theorem}

\begin{proof}
As in the proof of Proposition~\ref{prop:h1-P}, the Leray-Serre spectral
sequence implies that
$H^2(\U_P,V_\alpha) \cong [\Hcts^2(\u_P)\otimes V_\alpha]^R$.
It thus suffices to prove that the natural mapping
\begin{equation}
\label{eq:alg-galois}
H^2(\cG_P,V_\alpha)\to \Hcts^2(\G,V_\alpha)
\end{equation}
is injective. Let $[\gamma]$ be a class in the left hand 
side corresponding to a two-step extension 
$$
1 \to V_\alpha \to E_2 \to E_1 \to k \to 1.
$$
Since all $\cG_P$-modules are locally finite,
we may assume $E_2,E_1$ to be finite dimensional.
By exactness of the weight filtration and
the reductivity of $R$, we may assume that
$W_0E_1=E_1$, $W_{n(\alpha)}E_1=0$,
and that
$W_{-1}E_2=E_2$, $W_{n(\alpha)-1}E_2=0$.
We assume that $[\gamma]$ is a trivial
class as $\G$-module.
Then \cite[p.~575]{yoneda2} implies the
existence of a $\G$-module $E$ with 
$0 \to V_\alpha \to E \to E_1 \to 0$ and
$0 \to E_2 \to E \to k \to 0$.
Now the bootstrap property says that 
$E$ has property $P$, and hence is a weighted
$\G$-module with property $P$. Thus 
this is a $\cG_P$ module, which says
that $[\gamma]$ is the trivial class as an extension 
of $\cG_P$-modules, too.
Thus injectivity is proved.
\end{proof}

\subsection{The case where $\l$ is not contained in $S$}
\label{sec:crys}

We now show how to remove the hypothesis that $S$ contain all primes over
$\l$. Let $[\l]$ denote the set of primes of $F$ above the rational prime
$\l$. Assume that $S$ does not contain $[\l]$.  Then we define the category
$\MTMl(\OFS)$ of $\l$-adic mixed Tate modules over $\Spec\OFS$ to be the
full subcategory of $\MTMl(\O_{F, S\cup [\l]})$ consisting of the modules
$M$ with the property that for every prime $\p \in [\l] \backslash S$, the
representation $G_{F_\p} \to G_{F,S\cup[\l]} \to \Aut M$ is crystalline,
where $G_{F_\p}$ is the decomposition group of $G_F$ at the prime $\p$,
and $F_\p$ is the completion of $F$ at $\p$. (See \cite{FO,BK} for
crystalline representations.)

Let $P$ be the condition that a weighted $\GFSl$-module is crystalline
at every prime $\p$ in $[\l]$ outside $S$. 
It is known \cite{FO} that $P$ is closed under direct sums,
tensor products, taking subquotients and duals, and that the trivial
representation is  crystalline. We can therefore consider the weighted
completion of $\GFSl$ with respect to the cyclotomic character and
constrained by $P$ --- it is the tannakian fundamental group of
$\MTMl(\OFS)$. We denote this completion by $\AFS$, its prounipotent
radical by $\KFS$ and the Lie algebra of $\KFS$ by $\kFS$.

The group $\Hcts^1(\GFS,V)_P$ is the {\it finite part} ${\Hcts^1}_f(\GK, V)$
of $\Hcts^1(\GK,V)$, which is defined in \cite[(3.7.2)]{BK} and
\cite[p.~354]{BK} where it is denoted simply by $H^1_f(\GK,V)$.

One can show that the property $P$ has the bootstrap property
for $b=2$. The key point of the proof is that for a short exact
sequence
$$
0\to V_1 \to V_2 \to V_3 \to 0
$$
of crystalline $\l$-adic representations of $G_{F_\p}$, we have, by
\cite[Cor.~3.8.4]{BK}, a long exact sequence
\begin{eqnarray*}
0&\to& H^0(G_{F_\p},V_1)\to H^0(G_{F_\p},V_2)
\to H^0(G_{F_\p},V_3) \\ 
&&
\to {\Hcts^1}_f(G_{F_\p},V_1)
\to {\Hcts^1}_f(G_{F_\p},V_2)
\to {\Hcts^1}_f(G_{F_\p},V_3)
\to 0.
\end{eqnarray*}

\begin{proposition}
The graded Lie algebra $\Gr^W_\dot\kFS$ is free and
$$
\Hcts^1(\kFS) \cong \bigoplus_{n=1}^\infty
{\Hcts^1}_f(\GFS,\Ql(n))\otimes\Ql(-n)
\cong \bigoplus_{n=1}^\infty \Ql(-n)^{d_n},
$$
where $d_n$ is given as in Theorem~\ref{a_l}.
\end{proposition}

\begin{proof}
By Theorem~\ref{th:h2-P} we have an injection.
$$
\Hcts^2(\kFS) \hookrightarrow \bigoplus_{n=2}^\infty
{\Hcts^2}(\GFS,\Ql(n))\otimes\Ql(-n).
$$
Since the right hand side vanishes, $\Gr^W_\dot \kFS$ is free by
Lemma~\ref{free_crit}.

The computation of $H_1(\kFS)$ follows from Proposition~\ref{prop:h1-P}
and the facts (cf.\ Example~3.9 in \cite[p.~359]{BK}) that 
${\Hcts^1}_f(\GFS,\Ql(n))={\Hcts^1}(\GFS,\Ql(n))$ when $n\geq 2$ and that
there is an inclusion
$$
{\Hcts^1}_f(\GFSl,\Ql(1))\hookrightarrow {\Hcts^1}(\GFSl,\Ql(1))
$$
corresponding to $\O_{F,S}^\times \hookrightarrow \O_{F,S\cup[\l]}^\times$
via Kummer characters. Thus,
$$
\dim_\Ql{\Hcts^1}_f(\GFSl,\Ql(n)) = d_n
$$
where $d_n$ is defined as in Theorem~\ref{a_l}.
\end{proof}

For example, when $\OFS=\Spec \Z$,  $\Gr^W_\dot\kFS$ is a free Lie algebra
generated by $\sigma_3,\sigma_5,\sigma_7,\ldots$.

The formula for $\Ext^i_{\MTMl(\OFS)}(\Ql,\Ql(n))$
in Theorem~\ref{Th:MTMl} also holds for this case,
by replacing  $\Hcts^1(\GFS,\Ql(1))$ with ${\Hcts^1}_f(\GFS,\Ql(1))$.

\begin{remark}
Shiho \cite{shiho} has announced that his theory of the unipotent crystalline
fundamental group implies that the Lie algebra of $\Ql$-unipotent completion
of the fundamental group of a proper smooth curve is crystalline at $\p$ for
any prime above $\l$ if the curve has a good reduction at $\p$.

The corresponding result for incomplete varieties, including the projective
line minus three points, will follow from the syntomic conjecture for
incomplete varieties (cf.\ \cite{tsuji} for proper case) together with Shiho's
methods.

Instead, we established the vanishing of $\sigma_1$ using Lemma~\ref{lem:sigma1}. 
This and the above theorem show conversely that the $\Ql$-unipotent
completion of $\pi_1(\PminusQbar)$ is crystalline, since outer $\Gl$-action
factors through $\A_{\Q,\emptyset}$.
\end{remark}

\begin{theorem}
The outer action $\Gl \to \Out \cP$ of $\Gl$ on $\l$-adic unipotent
completion of $\pi_1(\PminusQbar)$ is crystalline at $\l$. \qed
\end{theorem}

\appendix
\section{Unipotent Completion}
\label{unipt_comp}

This appendix is a collection of results needed on unipotent completion.
Some of this material may be new, while other results are either folklore
or implicit in the literature. Unipotent completion is also known as
Malcev completion. One particularly useful approach to unipotent completion
is due to Quillen \cite[Appendix~A]{quillen}.

\subsection{Unipotent completion}
\label{basics}
Suppose that $\G$ is a group and that $k$ is a field of characteristic
zero. The {\it unipotent completion of $\G$ over $k$} consists of a 
prounipotent $k$-group $\G^\un_{/k}$ together with a homomorphism
$\theta : \G \to \G^\un_{/k}(k)$. These have the property that every
homomorphism $\G \to U(k)$ from $\G$ to the $k$-rational points of a
unipotent $k$-group $U$ factors uniquely through $\theta$.

Unipotent completion is easily seen to exist:
$$
\G^\un_{/k} = \limproj \rho U_\rho
$$ 
where $\rho$ ranges over all Zariski dense representations
$\rho : \G \to U_\rho(k)$ from $\G$ into a unipotent $k$-group.

The universal mapping property implies that unipotent completion is unique
up to a canonical isomorphism. The image of $\G$ in its unipotent completion
is Zariski dense, since the Zariski closure of the image has the requisite
universal mapping property.

Since all finite dimensional vector spaces are unipotent $k$-groups, there
is a canonical isomorphism
\begin{equation}
\label{h1_isom}
H_1(\G^\un_{/k}(k)) \cong H_1(\G,k):=\G^{ab}\otimes_\Z k
\end{equation}
whenever $H_1(\G,k)$ is finite dimensional.

Unipotent completion behaves well under field extension. The following
result is proved in \cite{hain:comp}.

\begin{proposition}
\label{ext_scalars}
Suppose that $\G$ is a group, that $k$ is a field of characteristic
zero and that $H_1(\G,k)$ is finite dimensional. If $K$ is an extension
field of $k$, then the natural homomorphism $\G \to \G^\un_{/k}(K)$
induces an isomorphism
$$
\G^\un_{/K} \stackrel{\sim}{\longrightarrow} \G^\un_{/k}\otimes_k K.
$$
In particular, if $\G$ is a finitely generated group, then there is a
natural isomorphism
$$
\G^\un_{/k} \stackrel{\sim}{\longrightarrow} \G^\un_{/\Q}\otimes_\Q k. \qed
$$
\end{proposition}

The lower central series
$$
\G = L^1 \G \supseteq L^2 \G \supseteq L^3 \G \supset \cdots
$$
of a group $\G$ is defined inductively by
$$
L^1 \G = \G \text{ and } L^{m+1}\G = [\G,L^m\G].
$$
Define a filtration $D^\dot_k$ of $\G$ by
$$
D^m_k \G =
\ker\big\{\G \to (\G^\un_{/k}/L^m\G^\un_{/k})(k)\big\}.
$$
This is a central filtration of $\G$, so it follows that
$L^m \G \subseteq D_k^m\G$ for all $m$.

\begin{theorem}
\label{lcs_gradeds}
If $H_1(\G,k)$ is finite dimensional, then the subgroup $D_k^m\G$ is the
inverse image of the set of torsion elements of $\G/L^m\G$ under the quotient
mapping $\G \to \G/L^m \G$. The central filtration $D^\dot_k$ of $\G$ is
the most rapidly descending central series with the property that each
$\G/D_k^m\G$ is a torsion free nilpotent group of length $< m$. If $\G$ is
finitely generated, then each $\Gr_{D_k}^m \G$ is a finitely generated abelian
group, and each $D^m_k\G /L^m \G$ is finite.
\end{theorem}

\begin{proof}
First recall the elementary fact that the set of torsion elements of
a nilpotent group $N$ forms a characteristic subgroup $T$ and that $N/T$ is
torsion free. This implies that the set $R^m\G$ of elements of $\G$ that
are torsion modulo $L^m \G$ is a characteristic subgroup, that $\G/R^m\G$
is a torsion free nilpotent group, and that $R^m\G/L^m \G$ is a torsion
group. It follows that the central series $R^\dot\G$ is the most rapidly
descending central series with each $\G/R^m\G$ torsion free.

Let $\cP = \G^\un_{/k}$.
Since $H_1(\G,k)$ is finite dimensional, each $\cP/L^m\cP$ is unipotent.
Since unipotent groups are torsion free, and since $\G/D^m_k\G$ is a
subgroup of $\cP(k)/L^m\cP(k)$, it follows that $D_k^m\G \supseteq R^m\G$.

It follows from Quillen's version \cite[Appendix A]{quillen} of Malcev's
Theorem \cite{malcev} that there is a unipotent $k$-group $U$ that
contains $\G/R^m\G$ as a Zariski dense subgroup. The density implies
that the length of $U$ is $<m$. The
homomorphism $\G \to U(\Ql)$ induces a homomorphism $\cP \to U$ which
factors through $\cP/L^m\cP$ as $U$ is unipotent of length $<m$. It 
follows that $R^m\G \supseteq D_k^m\G$, and that $R^m\G = D_k^m\G$ for
all $m$, which implies the result.
\end{proof}

\subsection{$\l$-adic unipotent completion}
\label{ladic_compln}

There is an obvious variant of unipotent completion for topological groups.
Suppose that $\G$ is a topological group and $k$ a topological field
with characteristic zero. We shall say that a homomorphism
$\G \to G(k)$ from $\G$ into the $k$-rational points of an algebraic
$k$-group $G$ is {\it continuous} if it is continuous with respect to
the topology on $G(k)$ induced from that of $k$. A homomorphism $\G \to G(k)$
from $\G$ to the $k$-points of a proalgebraic $k$-group $G$ is {\it continuous}
if it is the inverse limit of continuous homomorphisms $\G \to G_\alpha(k)$
from $\G$ to each of the finite dimensional quotients $G_\alpha$ of $G$.

The {\it continuous unipotent completion of $\G$} consists of a 
prounipotent $k$-group $\G^\un_{/k}$ together with a continuous
homomorphism $\theta : \G \to \G^\un_{/k}(k)$. 
These have the property
that every continuous homomorphism $\G \to U(k)$ from $\G$ to the
$k$-rational points of a unipotent $k$-group $U$ factors uniquely through
$\theta$.

The continuous unipotent completion of a topological group $\G$ is constructed
in much the same way as the standard unipotent completion given above. 
The universal mapping property of continuous unipotent completion  implies
that it is unique up to canonical isomorphism.

{}From now on, we assume that $\G$ is a profinite group and that $k=\Ql$
with the $\l$-adic topology. Since compact subgroups of $U(\Ql)$ are pro-$\l$
groups, the $\l$-adic unipotent completion of $\G$ factors through the
pro-$\l$ completion $\G^\prol$ of $\G$. Moreover, 
$$
\G^\un_\Ql \to (\G^\prol)^\un_\Ql 
$$
is an isomorphism since both have the same universal mapping property. Thus,
we may assume, without loss of generality, that $\G$ is a pro-$\l$ group.

The lower central series
$$
\G = L^1 \G \supseteq L^2 \G \supseteq L^3 \G \supset \cdots
$$
of a profinite group $\G$ is defined inductively by
$$
L^1 \G = \G \text{ and } L^{m+1}\G = \text{ closure of }[\G,L^m\G]
\text{ in } \G.
$$

Define a filtration $D_\l^\dot$ of $\G$ by
$$
D_\l^m \G = \ker\big\{\G \to (\G^\un_{/\Ql}/ L^m\G^\un_{/\Ql})(\Ql)\big\}.
$$
This is a central filtration of $\G$, so it follows that
$L^m \G \subseteq D_\l^m\G$ for all $m$.

\begin{theorem}
\label{ladic_lcs_gradeds} 
If $\G$ is  a topologically finitely generated pro-$\l$ group, then the
subgroup $D_\l^m\G$ is the inverse image of the set of torsion elements of
$\G/L^m\G$ under the quotient mapping $\G \to \G/L^m \G$. The central
filtration $D_\l^\dot$ of $\G$ is
the most rapidly descending central series with the property that each
$\G/D_\l^m\G$ is a torsion free nilpotent group of length $< m$. 
Moreover, each $\Gr_{D_\l}^m \G$ is a finitely generated
$\Zl$-module and each $D_\l^m\G /L^m \G$ is a finite $\l$-group.
\end{theorem}

\begin{corollary}
The $\l$-adic unipotent completions of $\G/L^m\G$ and $\G/D^m\G$ are both
$\G^\un_{/\Ql}/ L^m\G^\un_{/\Ql}$.
\end{corollary}

\begin{proof}
This follows by a proof that is essentially the same as that of
Theorem~\ref{lcs_gradeds}, taking the continuity of
$\G \to U(\Ql)$ into consideration in the current
setting. The following lemma is a key ingredient.
\end{proof}

\begin{lemma}
If $\G$ is a torsion-free nilpotent pro-$\l$ group, which is topologically
finitely generated, then $\G$ is continuously embeddable in the $\Ql$-rational
points of a unipotent group $U$ over $\Ql$.
\end{lemma}

\begin{proof}
The proof is similar to the proof of  Proposition~3.6 (a) in
\cite[Appendix A]{quillen}. Since $\G^{ab}$ is finitely generated, there is
a central filtration of $\G$ by closed normal subgroups, each of whose
graded quotients is either $\Zl$ or $\Z/\l$. We prove the lemma by induction
on the length of this filtration. There is an exact sequence
$$
1 \to \G_1 \to \G \to C \to 1,
$$
where $C$ is $\Zl$ or $\Z/\l$.
By induction, $\G_1$ is embedded in its $\l$-adic unipotent completion
$U_1$, which is algebraic by induction. If $C\cong \Zl$, then
$\G \cong C\ltimes \G_1$. The conjugacy action of $C$ on $\G_1$ lifts by
functoriality to a unipotent action on the Lie algebra of $U_1$, and thus
extends to a unipotent action of $\Ga$ on $U_1$.  It follows that
$C\ltimes \G_1$ can be continuously embedded in the unipotent group 
$\Ga \ltimes U_1$. If $C\cong \Z/\l$, then choose $x\in \G$ whose image
generates $C$.  Then $x^\l \in U_1(\Ql)$, and since $U_1(\Ql)$ is 
uniquely divisible, there is a unique $x'\in U_1(\Ql)$ with $x^\l = {x'}^\l$.
Then $\G$ can be imbedded in $U1(\Ql)$ by mapping $x$ to $x'$.
\end{proof}

\subsection{A comparison theorem}

There is a close relation between the unipotent completion of a finitely
generated group and the $\l$-adic unipotent completion of its profinite
and $\l$-adic completions.

\begin{theorem}\label{comp_thm}
If $\G$ is a finitely generated discrete group with pro-$\l$ completion
$\G^\prol$ and profinite completion $\Ghat$, then the three groups
$$
\G^\un_{/\Q}\otimes_\Q \Ql, \quad {\G^\prol}^\un_{/\Ql}
\text{ and } \Ghat^\un_{/\Ql}
$$
are all canonically isomorphic as prounipotent $\Ql$-groups.
\end{theorem}

The proof reduces to the following result:

\begin{lemma}
Suppose that $\G$ is a finitely generated group with pro-$\l$ completion
$\G^\prol$. If $U$ is a prounipotent group over $\Ql$, then every homomorphism
$\rho :\G \to U(\Ql)$ is continuous with respect to the pro-$\l$ topology on
$\G$, so that $\rho$ induces a continuous
homomorphism $\G^\prol \to U(\Ql)$. \qed
\end{lemma}

\begin{proof}
We may assume that $U$ is the upper triangular unipotent subgroup 
of $GL_N(\Ql)$ for some $N$.  We denote by $U(\l^m\Zl)$ the group of 
matrices whose $ij$th entry lies in $\l^{(i-j)m}\Zl$ when
$i>j$, is $1$ when $i=j$, and $0$ when $i<j$.
Since $\G$ is finitely generated, the image of $\rho$ is
contained inside $U(\l^m\Zl)$ for some $m\in \Z$. The
filtration
$$
\dots \supset U(\l^n\Zl) \supset U(\l^{n+1}\Zl)\supset \cdots 
$$
is a basic set of neighbourhoods of the identity in $U(\Ql)$; each quotient
is a finite group of $\l$-power order. Since $\G \subseteq U(\l^m\Zl)$,
the inverse image of each $U(\ell^n\Zl)$ is a finite index subgroup of $\G$
of $\l$-power order. The result follows.
\end{proof}

\begin{proof}[Proof of Theorem~\ref{comp_thm}]
It follows from the lemma that $\G \to \G^\un_{/\Q}(\Ql)$ is continuous and 
induces continuous homomorphisms
$$
\Ghat \to \G^\prol \to \G^\un_{/\Q}(\Ql).
$$
By the universal mapping property of $\l$-adic unipotent completion,
these induce homomorphisms
\begin{equation}\label{chain}
\Ghat^\un_{/\Ql} \to {\G^\prol}^\un_{/\Ql} \to \G^\un_{/\Q}\otimes_\Q \Ql.
\end{equation}
But the homomorphism $\G \to \Ghat \to \Ghat^\un_{/\Ql}(\Ql)$ induces a
homomorphism $\G^\un_{/\Q}\otimes_\Q \Ql \to \Ghat^\un_{/\Ql}$
whose composite with (\ref{chain}) is the canonical isomorphism of
Proposition~\ref{ext_scalars}. This completes the proof as the image
of $\G^\prol \to {\G^\prol}^\un_{/\Ql}(\Ql)$ is Zariski dense.
\end{proof}

\subsection{Automorphisms}
\label{autos}

Suppose that $\U$ is a prounipotent group over a field $k$ of
characteristic zero with Lie algebra $\u$. Standard Lie theory
implies that there is a natural isomorphism $H_1(\U) \cong H_1(\u)$.
Suppose that this vector space is finite dimensional over $k$. For each
field extension $K$ of $k$, define the group of continuous automorphisms
of $\U(K)$ by
$$
\Aut \U(K) = \limproj m \Aut(\u\otimes_k K/L^m\u\otimes_k K)
$$
where $L^m\u$ denotes the $m$th term of the lower central series of
$\u$.

\begin{proposition}
\label{algAutU}
If $H_1(\u)$ is finite dimensional over $k$, then there is a proalgebraic
group $\cAut(\U)$ defined over $k$ that represents the functor
$$
K \mapsto \Aut(\U(K))
$$
from field extensions of $k$ to groups. The kernel $\K$ of
the natural homomorphism $\cAut(\U) \to \Aut H_1(\U)$ is prounipotent.
The Lie algebra of $\cAut \U$ is isomorphic to $\Der \u$, the Lie algebra
of continuous derivations of $\u$.
\end{proposition}

\begin{proof}
Suppose that $\phi$ is an automorphism of $\u$. Since $H_1(\u)\cong\Gr_L^1\u$,
and since the bracket mapping
$$
\Gr_L^n \u \otimes \Gr_L^1 \u \to \Gr_L^{n+1} \u
$$
is surjective and commutes with $\phi$, we see that $\phi$ acts trivially
on $H_1(\u)$ if and only if it acts trivially on $\Gr_L^\dot\u$. It follows
that the kernel $\K$ of $\cAut \u \to \Aut H_1(\u)$
is a prounipotent group. The last statement follows from standard Lie
theory as $\Aut \U = \Aut \u$.
\end{proof}

Dividing out by the subgroup of inner automorphisms, we obtain:

\begin{corollary}
If $H_1(\u)$ is finite dimensional over $k$, then there is a 
proalgebraic group $\cOut(\U)$ defined over $k$ that represents the
functor
$$
K \mapsto \Out(\U(K))
$$
from field extensions of $k$ to groups. There is a canonical homomorphism
$\cOut \U \to \Aut H_1(\U)$ whose kernel is prounipotent. The Lie algebra
of $\cOut \U$ is naturally isomorphic to $\OutDer \u$, the Lie algebra of
continuous outer derivations of $\u$. \qed
\end{corollary}

\begin{corollary}
\label{cor:comparison}
For each finitely generated group $\G$, there are natural homomorphisms
$$
\Aut \G^\prol  \to \cAut (\G^\un_{/\Q}\otimes_\Q \Ql)
\to \cOut (\G^\un_{/\Q}\otimes_\Q \Ql).
$$
\end{corollary}

\begin{proof}
This follows from Theorem~\ref{comp_thm} as there is a natural
action of $\Aut \G^\prol$ on $\cAut ({\G^\prol}^\un_{/\Ql})$.
\end{proof}

{}From this we recover a result of Nakamura and Takao \cite{naka-taka}.

\begin{corollary}
Suppose that $k$ is a subfield of $\C$ and that $X$ is a variety over $k$.
Then for each $k$-rational point $x$ of $X$, there is a natural
homomorphism
$$
\Gal(\bar{k}/k) \to \cAut \pi_1(X(\C),x)^\un_{/\Q}\otimes\Ql. \qed
$$
\end{corollary}

\section{Continuous Cohomology of Galois Groups}
\label{gal_coho}

In this appendix we prove a result, stated below, that gives a computation
of the continuous cohomology of the Galois groups needed in this paper.
This is well known to the experts (cf.\ \cite{jannsen}), 
but we have included a proof for the convenience of readers.

\subsection{The result.}
Let $F$ be a number field, $\O_F$ its ring of integers, $\l$ a rational prime
number, $S$ a finite set of closed points of $\Spec \O_F$ containing all the
primes of $\O_F$ over $\l$. Let $\O_{F,S}$ be the ring of $S$-integers, so 
that $\Spec \O_{F,S}= \Spec \O_F - S$. Let $\GFS$ denote
$\pi_1(\Spec \O_{F,S}, \Spec \overline{F})$, the Galois group of the maximal
algebraic extension of $F$ unramified outside $S$. Let $r_1$ and $r_2$ be the
number of real and imaginary places of $F$, respectively.

The goal of this appendix is to prove the following result.

\begin{theorem}[Soul\'e]\label{Th:Galcoh}
With notation as above,
$$
\dim_\Ql
\Hcts^1(\GFS,\Ql(n)) =
\begin{cases}
r_1+r_2+\#S - 1 & \text{ if $n=1$,} \\
r_1+r_2 & \text{ if $n$ is odd and $\geq 3$,} \\
r_2 & \text{ if $n$ is even and positive.}
\end{cases}
$$
In addition, 
$\Hcts^2(\GFS,\Ql(n))$ vanishes for all $n\ge 2$.
\end{theorem}

The continuous cohomology $\Hcts^\dot$ is the one defined in
\cite[Sect.~2]{tate}. 
We reduce the proof to the computation of
the \'{e}tale cohomology of $\Spec \OFS$ by
the following proposition.

\begin{proposition}\label{prop:comparison}
For $i=1,2$ and all $n\in \Z$, there is a natural isomorphism
$$
\Hcts^i(\GFS,\Ql(n))\stackrel{\simeq}{\longrightarrow}
\Het^i(\Spec \O_{F,S},\Ql(n)). \qed
$$
\end{proposition}

In \'etale cohomology we have (cf.\ \cite{milne})
a canonical isomorphism
$$
\Het^i(\Spec \O_{F,S},\Ql(n))\cong
[\limproj{m}\Het^i(\Spec \O_{F,S},\Zl(n)/\l^m)]\otimes_\Zl \Ql.
$$
(For an abelian group $M$, let $M/\l^m = M/\l^m M$ and denote the
$\l^m$-torsion of $M$ by $M[\l^m]$.)
Since the canonical map
$$
\Hcts^i(\GFS,\Zl(n)/\l^m) \to \Het^i(\Spec \OFS,\Zl(n)/\l^m)
\ \ (i \geq 1)
$$
is isomorphic \cite[Proposition 2.9]{Milne2},
Proposition~\ref{prop:comparison} is a consequence
of the following lemma.

\begin{lemma}
\label{lem:limit}
When $i = 1,2$, there is a natural isomorphism
$$
\Hcts^i(\GFS,\Ql(n))\cong
[\limproj m \Hcts^i(\GFS,\Zl(n)/\l^m)]\otimes_\Zl \Ql. \qed
$$
\end{lemma}

Note that $\Hcts^\dot(G,\Zl(n)/{\l^m})$ is the usual cohomology of
profinite groups defined, for example, in \cite{serre2}.

\begin{proof}
The short exact sequence
$$
0 \to \Zl(n)\to \Ql(n) \to (\Ql/\Zl)(n) \to 0
$$
gives a long exact sequence \cite[p.~259]{tate}
\begin{multline*}
\cdots \to 
\Hcts^{i}(\GFS,\Zl(n)) \to
\Hcts^{i}(\GFS,\Ql(n)) \cr \to
\Hcts^{i}(\GFS,(\Ql/\Zl)(n)) \to
\Hcts^{i+1}(\GFS,\Zl(n)) \to
\cdots
\end{multline*}
Since $\GFS$ is compact, the proof of \cite[2.3]{tate}
shows that $\Hcts^{i}(\GFS,(\Ql/\Zl)(n))$ is torsion. So, after tensoring
with $\Ql$, we have an isomorphism
$$
\Hcts^i(\GFS,\Ql(n))\cong \Hcts^i(\GFS,\Zl(n))\otimes_\Zl \Ql.
$$

When $i=0,1$, $\Hcts^i(\GFS,\Zl(n)/\l^m)$ is finite. This is obvious when 
$i=0$. When $i=1$ it follows from  class field theory. Indeed, reduce
to the case where $\Zl(n)/\l^m$ has trivial Galois action by passing to
a finite extension. Then $\Hcts^1(\GFS,\Zl(n)/\l^m)$ is 
$\Homcts(\GFS,\Zl/\l^m)$, which is finite.
The corollary to \cite[Prop.~2.2]{tate} then implies that
$$
\Hcts^i(\GFS,\Zl(n))=
\limproj m \Hcts^i(\GFS,\Zl(n)/\l^m).
$$
Tensoring with $\Ql$ completes the proof.
\end{proof}

\subsection{Proof of Theorem~\ref{Th:Galcoh}}

First, we dispense with the
case where $n=1$ by showing that
\begin{equation}
\label{equality}
\dim_\Ql \Het^1(\Spec \O_{F,S},\Ql(1)) = r_1 + r_2 + \#S - 1. 
\end{equation}

Take the projective limit over $m$ of the Kummer sequence \cite[p.~128]{milne}
\begin{multline*}
0 \to \Het^0(\Spec \O_{F,S},\Gm)/\l^m \to 
\Het^1(\Spec \O_{F,S},\mu_{\l^m}) \cr \to
\Het^1(\Spec \O_{F,S}, \Gm)[\l^m] \to 0.
\end{multline*}
The limit of the center term is $\Het^1(\Spec \O_{F,S}, \Zl(1))$. Since
the ideal class group $\Het^1(\Spec \O_{F,S}, \Gm)$ is finite and the maps
in the inverse system are multiplication by powers of $\l$, the limit of the
right hand term vanishes. It follows that the limits of the first two terms
are isomorphic. Assertion (\ref{equality}) now follows as $\O_{F,S}^\times$
has rank $r_1 + r_2 + \#S - 1$ by the Dirichlet Unit Theorem.

To address the case where $n > 1$, we need the following result of Soul\'e
\cite[p.~376]{soule}.

\begin{theorem}[Soul\'e]
\label{th:soule}
If $F$ is a number field and $\l$ an odd prime number, then
$$
\dim_\Ql \Het^1(\Spec \O_F[\linv],\Ql(n)) =
\begin{cases}
r_1+r_2 & \text{ if $n$ is odd and $\geq 3$,} \\
r_2 & \text{ if $n$ is even and positive.}
\end{cases}
$$
In addition, $\Het^2(\Spec \O_F[\linv],\Ql(n))$ vanishes for all $n\ge 2$.
\end{theorem}

This theorem is also true when $\l=2$, following from
\cite{Rognes}. 
An explicit proof for the case of $F=\Q$
has been written down by Sharifi~\cite{sharifi}.

Theorem~\ref{Th:Galcoh} follows directly from this computation when $S$
is the set of primes lying over $\l$. For the general case, it is enough
to prove the following statement (cf.\ one may use 
\cite[Eqn.~13]{jannsen}).
 
\begin{lemma}
If $z$ is a closed point of $\Spec\O_{F,S}$ and $S'=S\cup \{z\}$, then
$$
\Het^i(\Spec \O_{F,S},\Ql(n))\cong \Het^i(\Spec \O_{F,S'},\Ql(n))
$$
for all $i\ge 0$ provided $n\neq 0,1$.
\end{lemma}

\begin{proof}
Let $X = \Spec \O_{F,S}$ and $U = \Spec \O_{F,S'}$. By
\cite[p.~92, Prop.~1.25]{milne} there is a long exact sequence
$$
\cdots \to 
H_z^i(X,\Ql(n))\to \Het^i(X,\Ql(n)) \to \Het^i(U,\Ql(n))
\to H_z^{i+1}(X,\Ql(n)) \to \cdots
$$
where $H_z^\dot$ denotes \'etale cohomology with support on $z$.
Thus it is enough to show that $H_z^i(X,\Ql(n))$ vanishes for all $i$ unless
$n=0$ or $1$.

Let $\X$ denote the henselization of $\Spec \O_F$ at $z$. Then by
\cite[p.~93 Cor.~1.28]{milne} we have
$$
H_z^i(X,\Ql(n))\cong H_z^i(\X,\Ql(n)),
$$
so it suffices to establish the vanishing of the right hand side.

Set $\U = \X-z$ and note that this is just $\Spec K(\X)$. By a second
application of the long exact sequence \cite[p.~92, Prop.~1.25]{milne}
applied to $(\X,\U)$, we see that if $\Het^i(\X,\Ql(n))$ and
$\Het^i(\U,\Ql(n))$ both vanish when $n\neq 0,1$, then $H_z^i(\X,\Ql(n))$
vanishes.

Firstly, by \cite[p.~224, Cor.~2.7]{milne}, we have
$$
\Het^i(\X,\Ql(n))=\Het^i(z,\Ql(n))=\Hcts^i(\hat{\Z},\Zl(n))\otimes \Ql.
$$
As in Example~\ref{zlx}, the right hand group is trivial unless $n=0$.

Set $\U^\ur = \Spec F(\X)$, where $F(\X)$ is the maximal algebraic
extension of $K(\X)$ unramified at $z$. To prove the second
vanishing we consider the Hochshild-Serre spectral sequence (cf.\
\cite[p.~105, Thm.~2.20]{milne})
$$
\Het^a(z, \Het^b(\U^\ur,\Ql(n)))\Rightarrow 
\Het^{a+b}(\U,\Ql(n)).
$$
According to \cite{serre2}, $\U^\ur$ has cohomological dimension at most 1,
so that
$$
\Het^b(\U^\ur,\Ql(n))=0 \text{ when } b > 1.
$$
Since the pro-$\l$ abelianization of the inertia group $\pi_1(\U^\ur)$ is
$\Z_l(1)$ as a Galois module, we have
$$
\Het^1(\U^\ur,\Ql(n)) = \Homcts(\pi_1(\U^\ur),\Ql(n)) \cong \Ql(n-1).
$$
The required vanishing follows by plugging this, and the fact that
$\Het^0(\U^\ur,\Ql(n))$ is $\Ql(n)$, into the spectral sequence.
\end{proof}

\end{document}